\documentclass[times]{nmeauth}

\usepackage{moreverb}

\usepackage[dvips,colorlinks,bookmarksopen,bookmarksnumbered,citecolor=red,urlcolor=red]{hyperref}
\usepackage{textcomp}
\usepackage{graphicx}
\usepackage{amssymb}
\usepackage{amsmath}
\usepackage[normalem]{ulem}
\usepackage{float}  %for figures location
\usepackage{threeparttable}

\newcommand\BibTeX{{\rmfamily B\kern-.05em \textsc{i\kern-.025em b}\kern-.08em
T\kern-.1667em\lower.7ex\hbox{E}\kern-.125emX}}

\def\volumeyear{2016}

\begin{document}

\runningheads{D.A. Bistrian and I.M. Navon}{Randomized Dynamic Mode Decomposition for NIROM}

\title{Randomized Dynamic Mode Decomposition for Non-Intrusive Reduced Order Modelling}%Reduced Order Models of Saint-Venant Systems

\author{D.A. Bistrian\affil{1}\corrauth\ and I.M. Navon\affil{2}}
\address{\affilnum{1}Department of Electrical Engineering and Industrial Informatics, Politehnica University of Timisoara, 331128 Hunedoara, Romania, \break
\affilnum{2}Department of Scientific Computing, Florida State University, Tallahassee, FL, 32306-4120, USA}

\corraddr{Revolutiei Str. Nr.5, 331128 Hunedoara, Romania. E-mail: diana.bistrian@upt.ro}

\begin{abstract}
This paper focuses on a new framework for reduced order modelling of non-intrusive data with application to $2D$ flows. To overcome the
shortcomings of intrusive model order reduction usually derived by combining the POD and the Galerkin projection methods, we developed a novel
technique based on Randomized Dynamic Mode Decomposition as a fast and accurate option in model order reduction of non-intrusive data originating
from Saint-Venant systems.  Combining efficiently the Randomized Dynamic Mode Decomposition algorithm with Radial Basis Function interpolation, we
produced an efficient tool in developing the linear model of a complex flow field described by non-intrusive (or experimental) data. The rank of
the reduced DMD model is given as the unique solution of a constrained optimization problem. We emphasize the excellent behavior of the
non-intrusive reduced order models by performing a qualitative analysis. In addition, we gain a significantly reduction of CPU time in computation
of the reduced order models (ROMs) for non-intrusive numerical data.

\end{abstract}


\maketitle

%\footnotetext[2]{Please ensure that you use the most up to date class file, available from the NME Home Page at\\
%\href{http://onlinelibrary.wiley.com/journal/10.1002/(ISSN)1097-0207}{\texttt{http://onlinelibrary.wiley.com/journal/10.1002/(ISSN)1097-0207}}}

\vspace{-6pt}

\section{Introduction, focus and motivation}
\vspace{-2pt} \label{intro}

The current line of this survey is motivated by the efficiency of reduced order modelling in different problems arising in hydrodynamics, where data
are collected in the aftermath of an experiment or are provided by measurement tools. We refer to this data as non-intrusive data. In the context of
application of mathematical and numerical techniques for modelling the non-intrusive massive data, the intent of this paper is to undertake an
efficient technique of model order reduction of shallow water systems. The challenge of this task is to build a linear dynamical system that models
the evolution of the strongly nonlinear flow  and to define a new mathematical and numerical methodology  for studying dominant and coherent
structures in the flow.

Among several model order reduction techniques, Proper Orthogonal Decomposition (POD) and Dynamic Mode Decomposition (DMD) represent modal
decomposition methods that are widely applied to study dynamics in different applications.
The application of POD is primarily limited to flows whose coherent structures can be hierarchically ranked in terms of their energy content. But
there are situations when the energy content is not a sufficient criterion to accurately describe the dynamical behavior of the flows. Instead, DMD
links the dominant flow features by a representation in the amplitudes-temporal dominant frequencies space.

Current literature has explored a broad variety of applications of reduced order modelling (ROM).
Recently, the POD approach has been incorporated for model order reduction purposes (MOR) by Chevreuil and Nouy \cite{Chevr2011}, Stefanescu and
Navon \cite{Stefanescu2013}, Dimitriu et al. \cite{dimistefanav2015}, Xiao et al. \cite{xiao2015}, Du et al. \cite{Du2013}, Fang et al.
\cite{Fang2013}. POD proved to be an effective technique embedded also in inverse problems, as it was demonstrated by the work of Cao et al.
\cite{cao2006, cao2007}, Chen et al. \cite{chen2011}, Stefanescu et al. \cite{stefasandunav2015}, thermal analysis (Bialecki et al.
\cite{Bialecki2004}), non-linear structural dynamics problems (Carlberg et al. \cite{Carlberg2011}) and aerodynamics (see for reference the recent
work of Semaan and his coworkers \cite{Semaan2016}). The method of Dynamic Mode Decomposition is rooted in the theory of Krylov subspaces
\cite{Chen2012, Rowley2009, Schmid2010} due to the original derivation of DMD as a variant of the Arnoldi algorithm \cite{Golub1996} and was first
introduced in 2008 by Schmid and Sesterhenn \cite{Schmid2008}. Only a year later, Rowley et al. \cite{Rowley2009} presented their theory of Koopman
spectral analysis, which has its inception back in 1931 \cite{Koopm1931}. In Dynamic Mode Decomposition, the modes are not orthogonal, but one
advantage of DMD compared to POD is that each DMD mode is associated with a pulsation, a growth rate and each mode has a single distinct frequency.
Owing to this feature, DMD method was used as a modal decomposition tool in non-linear dynamics (Schmid et al. \cite{Schmid2009}, Noack et al.
\cite{Noak2011}), in fluid mechanics (Bagheri \cite{Bagheri2013}, Rowley et al. \cite{Rowley2010}, Frederich and Luchtenburg \cite{Frederich2011},
Alekseev et al. \cite{aleksey2015}) and was recently introduced in turbulent flow problems (Seena and Sung \cite{Seena2011} and Hua et al.
\cite{Hua2016}). The theory of dynamic mode decomposition was applied for the purposes of model reduction within the efforts of Mezic
\cite{Mezic2005, Mezic2013} and in flow control problems by Bagheri \cite{Bagheri2012} and Brunton et al. \cite{Brunton2016}. For a complete
description of the utility of DMD vs. POD for model reduction in shallow water problems, the reader is referred to our previous paper (Bistrian and
Navon \cite{Bistrian2014}).

%

%%%%%%%%%%%%%%%%%%%%%%
The Saint-Venant system is employed in the present research to provide the experimental data. The Saint-Venant equations, named after the French
mathematician Adh\'{e}mar Jean Claude Barr\'{e} de Saint-Venant ($1797-1886$) (also called in literature the shallow-water equations (SWE)),
 are a set of hyperbolic partial differential equations that describe the flow below a pressure surface in a fluid. A
description of the Saint-Venant system has been presented by Vreugdenhil \cite{vre94}, as a result of depth-integration of the Navier–-Stokes
equations \cite{Galdi}. In literature, SWE are used in various forms to describe hydrological and geophysical fluid dynamics phenomena. The shallow
water magnetohydrodynamic (SWMHD) system has been devised by Gilman \cite{Gilman} to analyse the thin layer evolution of the solar tachocline.
Recently, a wave relaxation solver for SWMHD has been developed by Bouchut and Lh\'{e}brard \cite{Bouchut}. A collection of mathematical models in
coastal engineering based on SWE have been developed by Koutitus \cite{kou88}, while a reduced order modeling of the upper tropical Pacific ocean has
been presented by Cao et al. \cite{cao2006}.

The intrusive model order reduction is usually derived by combining the POD and the Galerkin projection methods \cite{Bistrian2014}. This approach
suffers from efficiency issues because the Galerkin projection is mathematically performed by laborious calculation and  requires stabilization
techniques in the process of numerical implementation, as it was argued in \cite{Amsallem2012, iliescu2014, Ballarin2015}.

The aim of this work is to circumvent these shortcomings by embedding the randomized dynamic mode decomposition as a fast and accurate option in
model order reduction of non-intrusive data originating from Saint-Venant systems. We propose in this work a novel  approach to derive a
non-intrusive reduced order linear model of the flow dynamics (denoted as NIROM) based on Dynamic Mode Decomposition of experimental data in
association with Radial Basis Function (RBF) multi-dimensional interpolation \cite{Golub1996}. Several key innovations are introduced in this paper
for the DMD-based model order reduction. The first one is represented by the randomization of the experimental data snapshots prior to singular value
decomposition of matrix data. Thus, we endow the DMD algorithm with a randomized SVD algorithm aiming to improve the accuracy of the reduced order
linear model and to reduce the CPU elapsed time. We gain a fast and accurate randomized DMD algorithm, exploiting an efficient low-rank DMD model of
input data. The rank of the reduced DMD model represents the unique solution of an optimization problem whose constraints are a sufficiently small
relative error of data reconstruction and a sufficiently high correlation coefficient
 between the experimental data and the DMD solution.

 We recall this procedure as adaptive randomized dynamic mode decomposition (ARDMD).

The first major advantage of the adaptive randomized DMD  proposed in this work is represented by the fact that the algorithm produces a reduced
order subspace of Ritz values, having the same dimension as the rank of randomized SVD function. In consequence, a further selection algorithm of the
Ritz values associated with their DMD modes is no longer needed. We employ in the flow reconstruction the smallest number of the DMD modes and their
amplitudes and Ritz values, respectively, leading to the minimum error of flow reconstruction, due to the adaptive feature of the proposed algorithm.
The second major improvement offered by the proposed randomized DMD can be found in the significantly reduction of CPU time for computation of
massive numerical data.

A randomized SVD algorithm was recently employed in conjunction with dynamic mode decomposition in \cite{Erichson2016} for processing of high
resolution videos in real-time. We believe that the present paper is the first work that shows the benefits of an efficient randomized dynamic mode
decomposition algorithm with application in fluid dynamics.

The reminder of the article is organized as follows. In Section 2 we recall the principles governing the Dynamic Mode Decomposition and we provide
the description of the DMD algorithm employed for decomposition of numerical data. In particular, we discuss the implementation of the randomized DMD
for optimal selection of the low order model rank. Section 3 is dedicated to theoretical considerations about multidimensional radial basis functions
interpolation. A detailed evaluation of the proposed numerical technique is presented in Section 4. Summary and conclusions are drawn in the final
section.

\vspace{-6pt}

\section{Adaptive randomized dynamic mode decomposition for non-intrusive data}
\vspace{-2pt} \label{section2}

\subsection{Koopman operator-the root of dynamic mode decomposition}

Being rooted in the work of French-born American mathematician Bernard Osgood Koopman \cite{Koopm1931}, the Koopman operator is applied to a
dynamical system evolving on a manifold $\mathbb{M}$ such that, for all ${v_k} \in \mathbb{M}$, ${v_{k + 1}} = f({v_k})$, and it maps any
scalar-valued function $g:\mathbb{M} \to \mathbb{R}$ into a new function ${\cal A}g$ given by
\begin{equation}
{\cal A}g\left( v \right) = g\left( {f\left( v \right)} \right).
\end{equation}

Has been demonstrated yet that spectral properties of the flow will be contained in the spectrum of the Koopman operator \cite{Bagheri2013} and even
when $f$  is finite-dimensional and nonlinear, the Koopman operator ${\cal A}$ is infinite-dimensional and linear \cite{Rowley2009, Mezic2013}. There
is a unique expansion that expands each snapshot in terms of vector coefficients ${\phi_j}$ which are called Koopman modes and mode amplitudes ${a
_j}\left( w \right)$, such that
iterates of ${v_0}$ are then given by
\begin{equation} \label{deco}
g\left( {{v_k}} \right) = \sum\limits_{j = 1}^\infty  {\lambda _j^k{a_j}} \left( {{v_0}} \right){\phi_j},\quad {\lambda _j} = {e^{{\sigma _j} + i{\omega _j}}},
\end{equation}
where ${\lambda _j}$ are called the Ritz eigenvalues of the modal decomposition, that are complex-valued flow structures associated with the growth
rate ${\sigma _j}$ and the frequency ${\omega _j}$.
Koopman modes represent spatial flow structures with time-periodic motion which are optimal in resolving oscillatory behavior. They have been
increasingly used  because they provide a powerful way of analysing nonlinear flow dynamics using linear techniques (see e.g. the work of Bagheri
\cite{Bagheri2013}, Mezic \cite{Mezic2013}, Rowley et al. \cite{Rowley2009}).
 The Koopman modes are extracted from the data snapshots and a unique frequency is associated to each mode. This is of major interest for fluid dynamics applications where phenomena occurring at different frequencies must be individualized.

A stable and consistent algorithm  proposed by Schmid \cite{Schmid2010}, referred in the literature as dynamic mode decomposition (DMD), can be used
for computing approximately a subset of the Koopman spectrum from the time series of snapshots of the flow. Thus DMD generalizes the global stability
modes and approximates the eigenvalues of the Koopman operator. The quantitative capabilities of DMD have already been well demonstrated in the
literature and several DMD procedures have been released  by the efforts of Rowley et al. \cite{Rowley2009}, Schmid \cite{Schmid2010}, Bagheri
\cite{Bagheri2013}, Mezic \cite{Mezic2013},  Belson et al. \cite{Belson2014}, Williams et al. \cite{Williams2015}.

Employing numerical simulations or experimental measurements techniques, different quantities associated with the flow are measured and collected as
observations at one or more time signals, called observables or non-intrusive data. It turns out (see the survey of Bagheri \cite{Bagheri2012}) that
monitoring an observable over a very long time interval allows the reconstruction of the flow phase space. Assuming that $\left\{
{{v_0},{v_1},...{v_N}} \right\}$ is a data sequence collected at a constant sampling time $\Delta t$, the DMD algorithm is based on the hypothesis
that a Koopman operator $\mathcal{A}$ exists, that steps forward in time the snapshots, such that
%such that
%\begin{equation}
%{v_{i + 1}} = \mathcal{A}{v_i},\quad i = 0,...,N - 1.
%\end{equation}
the snapshots data set
\begin{equation}
\left\{ {{v_0},\;{\cal A}{v_0},\;{{\cal A}^2}{v_0},.\;..,\;{{\cal A}^{N - 1}}{v_0}} \right\}
\end{equation}
corresponds to the $N^ {th}$ Krylov subspace generated by the Koopman operator from ${{v_0}}$.
The goal of DMD is to determine eigenvalues and eigenvectors  of the unknown matrix operator $\mathcal{A}$, thus a Galerkin projection of
$\mathcal{A}$ onto the subspace spanned by the snapshots is performed. For a sufficiently long sequence of the snapshots, we suppose that the last
snapshot ${{v_N}}$ can be written as a linear combination of previous $N-1$ vectors, such that
\begin{equation}
{v_N} = {c_0}{v_0} + {c_1}{v_1} + ... + {c_{N - 1}}{v_{N - 1}} + \mathcal{R},
\end{equation}
in which ${c_i},i = 0,...,N - 1$ are complex numbers and $\mathcal{R}$ is the residual vector. We assemble the following relations
\begin{equation}\label{kryl}
\mathcal{A}\left\{ {{v_0},{v_1},...{v_{N - 1}}} \right\} = \left\{ {{v_1},{v_2},...{v_N}} \right\} = \left\{ {{v_1},{v_2},...,V_0^{N - 1}\textit{c}} \right\} + {\mathop{\rm \mathcal{R}e}\nolimits} _{N - 1}^T,
\end{equation}
where $V_0^{N - 1} = \left( {\begin{array}{*{20}{c}} {{v_0}}&{{v_1}}&{...}&{{v_{N - 1}}}
\end{array}} \right)$, ${\textit{c}^T} = \left( {\begin{array}{*{20}{c}}{{c_0}}&{{c_1}}&{...}&{{c_{N - 1}}}\end{array}} \right)$ is the unknown complex column vector and $e_j^T$ represents the $j^ {th}$ Euclidean unitary vector of length $N−-1$.

In matrix notation form, Eq. (\ref{kryl}) reads
\begin{equation}\label{rel}
\mathcal{A}V_0^{N - 1}  = V_0^{N - 1}\mathcal{C} + {\mathop{\rm \mathcal{R}e}\nolimits} _{N - 1}^T,\quad \mathcal{C} = \left( {\begin{array}{*{20}{c}}
0&{...}&0&{{c_0}}\\
1&{}&0&{{c_1}}\\
 \vdots & \vdots & \vdots & \vdots \\
0& \ldots &1&{{c_{N - 1}}}
\end{array}} \right),
\end{equation}
where $\mathcal{C}$ is the companion matrix.

A direct consequence of (\ref{rel}) is that decreasing the residual increases overall convergence and therefore the eigenvalues of the companion
matrix $\mathcal{C}$ will converge toward the eigenvalues of the Koopman operator $\mathcal{A}$.

The representation of data in terms of Dynamic Mode Decomposition is given by the linear model
\begin{equation}\label{dmdrep}
{v_{DMD}}^t\left( \textbf{x} \right) = \sum\limits_{j = 1}^k {{a_j}{\phi _j}\left( \textbf{x} \right)\lambda _j^{t - 1}} ,\quad {\lambda _j} = {e^{\left( {{\sigma _j} + i{\omega _j}} \right)\Delta t}},\quad t = t_1,...,t_N,
\end{equation}
where the right eigenvectors ${\phi _j} \in \mathbb{C}$ are dynamic (Koopman) modes, the eigenvalues  ${\lambda _j}$ are called Ritz values
\cite{Chopr2000} and coefficients ${a_j} \in \mathbb{C}$ are denoted as amplitudes or Koopman eigenfunctions. Each Ritz value ${\lambda _j}$ is
associated with the growth rate ${\sigma _j} = \frac{{\log \left( {\left| {{\lambda _j}} \right|} \right)}}{{\Delta t}}$ and the frequency ${\omega
_j} = \frac{{\arg \left( {\left| {{\lambda _j}} \right|} \right)}}{{\Delta t}}$ and  $k$ represents the number of DMD modes involved in
reconstruction.

Since it was first introduced in 2008 by Schmid and Sesterhenn \cite{Schmid2008}, a considerable amount of work has focused on understanding and
improving the method of dynamic mode decomposition. Chen et al. \cite{Chen2012} introduced an optimized DMD, which tailors the decomposition to an
optimal number of modes. This method minimizes the total residual over all data vectors and uses simulated annealing and quasi-Newton minimization
iterative methods for selecting the optimal frequencies. A recursive dynamic mode decomposition was developed by Noack et al. \cite{noack2015} with
application to a transient cylinder wake. Multi-Resolution DMD was released by Kutz et al. \cite{kutz2015} for extracting DMD  eigenvalues and modes
from data sets containing multiple timescales. Efficient post processing procedures for selection of the most influential DMD modes and eigenvalues
were presented in our previous papers \cite{Bistrian2014, bisnavon2016}.

So far we have noticed two directions in developing the algorithms for dynamic mode decomposition. The straight-forward approach is seeking the
companion matrix $\mathcal{C}$ from (\ref{rel}) that helps to construct in a least squares sense the final data vector as a linear combination of all
previous data vectors \cite{Rowley2009, Fiedler2003, Rowley2010}. Because this version may be ill-conditioned in practice, Schmid \cite{Schmid2010}
recommends an alternate algorithm, based on Singular Value Decomposition (SVD) \cite{Golub1996} of snapshot matrix, upon which the work within this
article is based.

\subsection{Adaptive randomized DMD algorithm}

To derive the improved algorithm proposed here, we proceed by collecting data ${v_i}\left( {t,\textbf{x}} \right) = v\left( {{t_i},\textbf{x}}
\right),\;{t_i} = i\Delta t,\;i = 0,...,N$, $\textbf{x}$ representing the spatial coordinates whether Cartesian or Cylindrical and form the snapshot
matrix $V = \left[ {\begin{array}{*{20}{c}} {{v_0}}&{{v_1}}&{...}&{{v_N}}
\end{array}} \right]$.

We arrange the snapshot matrix into two matrices. A matrix $V_0^{N - 1}$ is formed with the first $N$ columns and the matrix $V_1^N$ contains the
last $N$ columns of $V$: $V_0^{N - 1} = \left[ {\begin{array}{*{20}{c}} {{v_0}}&{{v_1}}&{...}&{{v_{N - 1}}}
\end{array}} \right]$, $V_1^N = \left[ {\begin{array}{*{20}{c}}
{{v_1}}&{{v_2}}&{...}&{{v_N}}
\end{array}} \right]$.

Expressing $V_1^N$ as a linear combination of the independent sequence $V_0^{N - 1}$ yields:
\[V_1^N = \mathcal{A}V_0^{N - 1} = V_0^{N - 1}\mathcal{S} + R,\]
where $R$ is the residual matrix and $\mathcal{S}$ approximates the eigenvalues of $\mathcal{A}$ when ${\left\| R \right\|_2} \to 0$. The objective
at this step is to solve the minimization problem
\begin{equation}\label{minpro}
\mathop {Minimize}\limits_\mathcal{S} \;R = {\left\| {V_1^N - V_0^{N - 1}\mathcal{S}} \right\|_2}.
\end{equation}

An estimate can be computed by multiplying $V_{1}^{{{N}}}$ by the Moore-Penrose pseudoinverse of $V_{0}^{{{N}}-1}$:
\begin{equation}
\mathcal{S}={{\left( V_{0}^{{{N}}-1} \right)}^{+}}V_{1}^{{{N}}}=W{{\Sigma }^{+}}{{U}^{H}}V_{1}^{{{N}}}=X\Lambda {{X}^{-1}},
\end{equation}
where $X$ and $\Lambda$ represent the eigenvectors, respectively the eigenvalues of $\mathcal{S}$, and ${{\Sigma }^{+}}$ is computed according to
Moore-Penrose pseudoinverse definition of Golub and van Loan \cite{Golub1996}:
\begin{equation}
{{\Sigma }^{+}}=diag\left( \frac{1}{{{\sigma }_{1}}},\cdots ,\frac{1}{{{\sigma }_{r}}},0\cdots ,0 \right),\quad r=rank\left( V_{0}^{{{N}}-1} \right).
\end{equation}

It can be seen that the SVD plays a central role in computing the DMD. Therefore, this approach for computing the operator $\mathcal{S}$ previously
employed in \cite{Bistrian2014} might not be feasible when dealing with high dimensional non-intrusive data.  It is more desirable to reduce the
problem dimension to avoid a computationally expensive SVD. Several key innovations are introduced in this paper for the DMD-based model order
reduction.

First one is represented by the randomization of the experimental data matrix $V_{0}^{{{N}}-1}$ prior to singular value decomposition. Thus, we endow
the DMD algorithm with a randomized SVD algorithm aiming to improve the accuracy of the reduced order linear model and to reduce the CPU elapsed
time. We gain a fast and accurate randomized DMD algorithm, exploiting an efficient low-rank DMD model of input data.

\textbf{Algorithm 1} describes the procedure for computing the randomized SVD and it is adapted after Halko et al. \cite{halko2011}.

\noindent\rule{12.5cm}{0.4pt}

\textbf{Algorithm 1:} Randomized SVD algorithm (RSVD)

\noindent\rule{12.5cm}{0.4pt}

\textbf{Initial data:} $V_0^{N - 1} \in {\mathbb{R}^{m \times n}}$, $m \ge n$, integer target rank $k \ge 2$ and $k < n$.

\begin{enumerate}

\item[1.] Generate random matrix $M = rand\left( {n,r} \right)$, $r = \min \left( {n,2k} \right)$.

\item[2.] Multiplication of snapshot matrix with random matrix $Q = V_0^{N - 1}M$.

\item[3.] Orthogonalization $Q \leftarrow orth\left( Q \right)$.

\item[4.] Projection of snapshot matrix $V = {Q^H}V_0^{N - 1}$, where $H$ denotes the conjugate transpose.

\item[5.] Produce the economy size singular value decomposition $\left[ {{Q_1},\Sigma ,W} \right] = SVD\left( V \right)$.

\item[6.] Compute the right singular vectors $U = Q{Q_1}$.

\textbf{Output:}  Procedure returns $U \in {\mathbb{R}^{m \times k}}$, $\Sigma  \in {\mathbb{R}^{k \times k}}$, $W \in {\mathbb{R}^{n \times k}}$.

\end{enumerate}

We define the relative error of the low-rank model as the ${L_2}$-norm of the difference between the experimental variables  and approximate DMD
solutions over the exact one, that is,
\begin{equation}\label{eror}
E{r_{DMD}} = \frac{{{{\left\| {v\left( {\textbf{x}} \right) - {v_{DMD}}\left( {\textbf{x}} \right)} \right\|}_2}}}{{{{\left\| {v\left( {\textbf{x}} \right)} \right\|}_2}}},
\end{equation}
where $v\left( {\textbf{x}} \right)$ represent the experimental data and ${v_{DMD}}\left( {\textbf{x}} \right)$ represent the low-rank DMD
approximation.

The correlation coefficient defined below is used as additional metric to validate the quality of the low-rank DMD model:
\begin{equation}\label{corelation}
{C_{DMD}} = \frac{{{{\left( {{{\left\| {v\left( \textbf{x} \right) \cdot {v_{DMD}}\left( \textbf{x} \right)} \right\|}_F}} \right)}^2}}}{{{{\left\| {v{{\left( \textbf{x} \right)}^H} \cdot v\left( \textbf{x} \right)} \right\|}_F}{{\left\| {{v_{DMD}}{{\left( \textbf{x} \right)}^H} \cdot {v_{DMD}}\left( x \right)} \right\|}_F}}},
\end{equation}
where $v\left( \textbf{x} \right)$ means the experimental data, ${{v_{DMD}}\left( x \right)}$  represent the computed solution by means of the
reduced order DMD model, $\left(  \cdot  \right)$ represents the Hermitian inner product and $H$ denotes the conjugate transpose. We denote by
${\left\| {\, \cdot \,} \right\|_F}$ the Frobenius matrix norm in the sense that for any matrix $A \in {\mathbb{C}_{m \times n}}$ having singular
values ${\sigma _1},...,{\sigma _n}$ and SVD of the form $A = U\Sigma {V^H}$, then
\begin{equation}
{\left\| A \right\|_F} = {\left\| {{U^H}AV} \right\|_F} = {\left\| \Sigma  \right\|_F} = \sqrt {{\sigma _1}^2 + ...+{\sigma _n}^2}.
\end{equation}

The rank of the reduced DMD model is given such that the relative error of data reconstruction becomes sufficiently small and the correlation
coefficient is sufficiently high. We recall this procedure as \textit{adaptive randomized dynamic mode decomposition}. Determination of the optimal
rank of the reduced DMD model then amounts to finding the solution to the following optimization problem:
\begin{equation}\label{optimprob}
\left\{ {\begin{array}{*{20}{l}}
{\mathop {Find}\limits_{k \in \mathbb{N},{\kern 1pt} k \ge 2} \;\;{v_{DMD}}^t\left( {\bf{x}} \right) = \sum\limits_{j = 1}^k {{a_j}{\phi _j}\left( {\bf{x}} \right)\lambda _j^{t - 1}} ,}\\
{Subject\;to\;\;k = \arg \min E{r_{DMD}}\;and\;k = \arg \max {C_{DMD}}.}
\end{array}} \right.
\end{equation}

The adaptive randomized DMD algorithm (\textbf{Algorithm 2}) that we applied in the forthcoming section for data originating from the Saint-Venant
system proceeds as follows:

\noindent\rule{12.5cm}{0.4pt}

\textbf{Algorithm 2:} Adaptive randomized DMD algorithm (ARDMD)

\noindent\rule{12.5cm}{0.4pt}

\textbf{Initial data:} $V_0^{N - 1} \in {\mathbb{R}^{m \times n}}$, $V_1^{N} \in {\mathbb{R}^{m \times n}}$, $m \ge n$, integer target rank $k \ge 2$
and $k < n$.

\begin{enumerate}

\item[1.] Variate $k$.

\item[2.] Produce the randomized singular value decomposition: $\left[ {U,\Sigma ,W} \right] = RSVD\left( V_0^{N - 1},k \right)$, where $U$
    contains the proper orthogonal modes of $V_0^{N - 1}$ and $\Sigma $ contains the singular values.

\item[3.] Solve the minimization problem (\ref{minpro}): $S = {U^H}\left( {V_1^NW{\Sigma ^{ - 1}}} \right)$.

For the reader information we will detail it in the following. Relations $\mathcal{A}V_0^{N - 1} = V_1^N = V_0^{N - 1}S + R, {\left\| R \right\|_2}
\to 0$ and $V_0^{N - 1} = U\Sigma {W^H}$ yield:
      \[\begin{array}{l}
\mathcal{A}U\Sigma {W^H} = V_1^N = U\Sigma {W^H}S\quad \\
\quad \quad \quad \quad  \Rightarrow \quad {U^H}\mathcal{A}U\Sigma {W^H} = {U^H}U\Sigma {W^H}S\quad \\
\quad \quad \quad \quad  \Rightarrow \quad \;S = {U^H}\mathcal{A}U.
\end{array}\]
From $\mathcal{A}U\Sigma {W^H} = V_1^N$ it follows that $\mathcal{A}U = V_1^NW{\Sigma ^{ - 1}}$ and hence $S = {U^H}\left( {V_1^NW{\Sigma ^{ - 1}}}
\right)$.

\item[4.]	Compute dynamic modes solving the eigenvalue problem $SX = X\Lambda $ and obtain dynamic modes as $\Phi  = UX$. The diagonal entries of
    $\Lambda$ represent the eigenvalues $\lambda$.

\item[5.]  Project dynamic modes onto the first snapshot to calculate the vector containing dynamic modes amplitudes $Ampl = \left( {{a_j}}
    \right)_{j = 1}^{rank\left( \Lambda  \right)}$.

\item[6.] The DMD model of rank $k$ is given by the product
\begin{equation}
V_{DMD} = \Phi  \cdot diag\left( Ampl \right) \cdot Van,
\end{equation}
where the Vandermonde matrix is
\[Van = \left( {\begin{array}{*{20}{c}}
1&{\lambda _1^1}&{\lambda _1^2}& \ldots &{\lambda _1^{N - 2}}\\
1&{\lambda _2^1}&{\lambda _2^2}& \ldots &{\lambda _2^{N - 2}}\\
1& \vdots & \vdots & \vdots & \vdots \\
 \ldots & \ldots & \ldots & \ldots & \ldots \\
1&{\lambda _k^1}&{\lambda _k^2}& \ldots &{\lambda _k^{N - 2}}
\end{array}} \right).\]

\item[7.] Solve the optimization problem (\ref{optimprob}) and obtain the optimal low rank $k$ and associated $V_{DMD}$.

\textbf{Output:} $k$, $V_{DMD}$.
\end{enumerate}

The first major advantage of the adaptive randomized DMD  proposed in this paper is represented by the fact that \textbf{Algorithm 2} produces a
reduced order subspace of Ritz values, having the same dimension as the rank of RSVD function. In consequence, after solving the optimization problem
(\ref{optimprob}), an additional selection algorithm of the Ritz values associated with their DMD modes is no longer needed. We employ in the flow
reconstruction the most influential DMD modes associated with their amplitudes and Ritz values, respectively, leading to the minimum error of flow
reconstruction, due to the adaptive feature of the proposed algorithm.

The second major improvement offered by the proposed randomized DMD can be found in the significantly reduction of CPU time for computation of
massive numerical data, as we will detail in the section dedicated to numerical results.

\section{Multidimensional radial basis function interpolation}\label{RBF}

\subsection{Model-order reduction using projection}

So far, the model order reduction practitioners applied the intrusive model order reduction having modal decomposition (POD or DMD) and the Galerkin
projection compound.

%We consider a bounded open domain $\Omega  \subset \mathbb{R}^3$ and let ${L^2}\left( \Omega  \right)$ be the Hilbert space of square integrable vector functions over $\Omega$, associated with the energy norm ${\left\| w \right\|_{{L^2}}} = \left( {w,w} \right)_{{L^2}}^{1/2}$ and the standard inner product ${\left( {v,w} \right)_{{L^2}}} = \int\limits_\Omega  {v \cdot w\,dz} $.
%Let ${H_\nabla }$ be the Hilbert space of divergence free functions given by
%%
%\begin{equation}
%{H_\nabla } = \left\{ {w \in {L^2}\left( \Omega  \right)\;\left| {\nabla  \cdot w = 0\;in\;\Omega ,\;w \cdot \mathop n\limits^ \to   = 0\;on\;\partial \Omega } \right.} \right\},
%\end{equation}
%%
%where ${\mathop n\limits^ \to  }$ is the outward normal to the boundary. We define ${H^d}\left( \Omega  \right) \subset {L^2}\left( \Omega  \right)$ to be the Hilbert space of functions $w$ with $d$ distributional derivatives $\nabla {w_i}$, $1 \le i \le d$, that are all square integrable. Let $\mathbb{V}$ be the Hilbert space
%%
%\begin{equation}
%\mathbb{V} = \left\{ {w \in {H_\nabla }\;\left| {w \in {H^1}\left( \Omega  \right),\;w = 0,\;\frac{{\partial w}}{{\partial \mathop n\limits^ \to  }} = 0,\;on\;\partial \Omega } \right.} \right\},
%\end{equation}
%%
%with norm ${\left\| w \right\|_\mathbb{V}} = \left( {w,w} \right)_\mathbb{V}^{1/2}$ and the inner product ${\left( {v,w} \right)_\mathbb{V}} = \sum\limits_{i = 1}^d {\left( {\nabla {v_i},\nabla {w_i}} \right)} $,

In the Cartesian or Cylindrical coordinates formulation, we suppose there exists a time dependent flow $v = \left( {\bf{x},t} \right) \in \mathbb{V}$
and a given initial flow $v\left( {\bf{x},0} \right)$, that are solutions of the coupled system of nonlinear ordinary differential equations
\begin{equation}\label{compact}
\left\{ \begin{array}{l}
\frac{{\partial v}}{{\partial t}}\left( {\bf{x},t} \right) = f\left( {t,v\left( {\bf{x},t} \right)} \right),\\
v\left( {\bf{x},{t_0}} \right) = {v_0}\left( {\bf{x}} \right),
\end{array} \right.
\end{equation}
obtained from the spatial discretization of evolution equations in continuous space.

To perform the intrusive model order reduction, we start by replacing the field $v$ with $v_{DMD}$ in (\ref{compact}) and project the resulting
equations onto the subspace ${X^{DMD}} = span\left\{ {{\phi _1}( \cdot ),{\phi _2}( \cdot ),...,{\phi _{k}}( \cdot )} \right\}$ spanned by the DMD
basis to compute the following inner products:
\begin{equation}\label{dmdrom}
\left\langle {{\phi _i}\left(  \cdot  \right),\sum\limits_{j = 1}^{k} {{\lambda _j}{\phi _j}\left(  \cdot  \right){{\dot a}_j}\left( t \right)} } \right\rangle  = \left\langle {{\phi _i}\left(  \cdot  \right),f\left( {t,\sum\limits_{j = 1}^{k} {{\lambda _j}{\phi _j}\left(  \cdot  \right){a_j}\left( t \right)} } \right)} \right\rangle  ,
\end{equation}
\begin{equation}
\left\langle {{\phi _i}\left(  \cdot  \right),\sum\limits_{j = 1}^{k} {{\lambda _j}{\phi _j}\left(  \cdot  \right){{\dot a}_j}\left( {{t_0}} \right)} } \right\rangle  = \left\langle {{\phi _i}\left(  \cdot  \right),{w_0}} \right\rangle ,\; \mbox{for} \; i = 1,...,k,
\end{equation}
where $\left\langle {f,g} \right\rangle  = \int_\Omega  {fg\,d\Omega } $.

The Galerkin projection gives the DMD-ROM, i.e., a dynamical system for temporal coefficients ${\left\{ {{a_j}\left( t \right)} \right\}_{j =
1,...,k}}$:
\begin{equation}\label{dmdcoef}
{\dot a_i}\left( t \right) = \left\langle {{\phi _i}\left(  \cdot  \right),f\left( {t,\sum\limits_{j = 1}^{k} {{\lambda _j}{\phi _j}\left(  \cdot  \right){a_j}\left( t \right)} } \right)} \right\rangle ,
\end{equation}
with the initial condition
\begin{equation}
{a_i}\left( {{t_0}} \right) = \left\langle {{\phi _i}\left(  \cdot  \right),{v_0}} \right\rangle ,\; \mbox{for} \; i = 1,...,k.
\end{equation}

The resulting autonomous system has linear and quadratic terms parameterized by ${{c_{im}}}$, ${{c_{imn}}}$, respectively:
\begin{equation}
{\dot a_i}\left( t \right) = \sum\limits_{m = 1}^{k} {\sum\limits_{n = 1}^{k} {{c_{imn}}} } {a_m}\left( t \right){a_n}\left( t \right) + \sum\limits_{m = 1}^{k} {{c_{im}}} {a_m}\left( t \right),\quad i = 1,...,k.
\end{equation}

DMD in combination with the Galerkin projection method is an effective method for deriving a reduced order model. For a detailed description of this
method, the reader is invited to consult our previous paper \cite{Bistrian2014}. By projecting the full dynamical system onto a reduced space which
is constructed based on the optimal DMD basis functions, the computational efficiency can be enhanced by several orders of magnitude. However, this
approach presents several shortcomings. As in the case of POD-ROMs \cite{Amsallem2012, iliescu2014, Ballarin2015}, this method implies analytical
calculations  and it remains dependent on the governing equations of the full physical system, therefore is not applicable in case of non-intrusive
(or experimental) data.

Recently, a diversity of non-intrusive methods have been introduced into ROMs, associated so far with proper orthogonal decomposition, like: Smolyak
sparse grid method \cite{XiaoLin2016}, radial basis functions interpolation \cite{Walton2013, BisResiga2016, XiaoYahng2016} or the method of
generalized moving least squares \cite{Deh2016}.

 In this paper we propose a novel approach to derive a reduced order model for non-intrusive data. In the offline stage of the proposed technique, the method of randomized dynamic mode decomposition finds the subspace  ${X^{DMD}} = span\left\{ {{\phi _1}( \cdot ),{\phi _2}( \cdot ),...,{\phi _{k}}( \cdot )} \right\}$ spanned by the sequence of the most efficient DMD modes. In the online stage, we involve the effectual multidimensional radial basis function interpolation who elegantly approximate the values of temporal coefficients for new time instances.

 \subsection{Model-order reduction by radial basis function interpolation}

Since it was introduced by Rolland Hardy in $1970$ \cite{Hardy1970} for applications in cartography, radial basis function (RBF) method  has
undergone a rapid progress  as an active tool of mathematical interpolation of scattered data in many application domains like domain decomposition
\cite{Beatson2000}, unsteady fluid flows modelling \cite{Walton2013, xiao2015} or image processing \cite{Carr2001}. Investigations upon accuracy and
stability of RBF based interpolation may be found in Fornberg and Wright \cite{Fornberg2004}, Fasshauer \cite{Fasshauer2007} and  Chenoweth
\cite{Chenoweth2009}.

The ARDMD algorithm previously described allows the identification of a reduced order model of form
\begin{equation}\label{dmdrep}
{v_{DMD}}^t\left( \textbf{x} \right) = \sum\limits_{j = 1}^k {{a_j}{\phi _j}\left( \textbf{x} \right)\lambda _j^{t - 1}} ,\quad {\lambda _j} = {e^{\left( {{\sigma _j} + i{\omega _j}} \right)\Delta t}},\quad t = t_1,...,t_N,
\end{equation}
in which ${\phi _j} \in \mathbb{C}$ represents dynamic DMD modes, ${\lambda _j}$ are the Ritz values, ${a_j} \in \mathbb{C}$ represent the modal
amplitudes and $k$ is the truncation order.

In the following, we employ the method of RBF interpolation, as a generalization of Hardy's multiquadric and inverse multiquadric method
\cite{Hardy1990}, for numerical interpolation of the model coefficients ${b_j}^t = {a_j}\lambda _j^{t - 1}$ for $t \in \left[ {{t_1},{t_N}} \right]$.

Considering the determined coefficients as a set of distinct nodes $\left\{ {{\boldsymbol x_i}} \right\}_{i = 1}^{k \times N} \subset {\mathbb{R}^2}$
and a set of function values $\left\{ {{f_i}} \right\}_{i = 1}^{k \times N} \subset \mathbb{R}$, the problem reduces to find an interpolant
$s:{\mathbb{R}^2} \to \mathbb{R}$ such that
\begin{equation}\label{interpcond}
s\left( {{\boldsymbol x_i}} \right) = {f_i}\quad for\;i = 1,...,k \times N,
\end{equation}
where $N$ is the number of time instances for which experimentally measured data are available and $k$ is the number of retained DMD modes. Note that
we use the notation ${f_i} = {b_j}^t$, $j = 1,...,k$, $t  = t_1,...,t_N$ for scattered points values and $\boldsymbol x = \left( {x,y} \right) \in
\left\{ {1,...,k} \right\} \times \left[ {{t_1},{t_N}} \right]$ for scattered points coordinates.

Considering $B{L_2}\left( {{\mathbb{R}^2}} \right)$ the Beppo-Levi space \cite{Aikaw1988} of distributions on ${\mathbb{R}^2}$ with square integrable
second derivatives, equipped with the rotation invariant semi-norm
\begin{equation}\label{seminorm}
 {\left\| s\left( {\boldsymbol{x}} \right) \right\|^2} = \int_{{\mathbb{R}^2}} {{{\left( {\frac{{{\partial ^2}s\left( \boldsymbol x \right)}}{{\partial {x^2}}}} \right)}^2} + } {\left( {\frac{{{\partial ^2}s\left( \boldsymbol x \right)}}{{\partial {y^2}}}} \right)^2} + 2{\left( {\frac{{{\partial ^2}s\left( \boldsymbol x \right)}}{{\partial x\partial y}}} \right)^2}d\boldsymbol x
\end{equation}
we seek the smoothest interpolant surface in the affine space
\begin{equation}\label{space}
S_{BL} = \left\{ {s \in B{L_2}\left( {{\mathbb{R}^2}} \right)\;\left| {\;s\left( {{\boldsymbol x_i}} \right) = {f_i},\quad i = 1,...,k \times N} \right.} \right\},
\end{equation}
i.e.,
\begin{equation}\label{argmin}
\widetilde s\left( {\boldsymbol{x}} \right) = \arg \;\mathop {\min }\limits_{s \in {S_{BL}}} {\left\| s\left( {\boldsymbol{x}} \right) \right\|^2}.
\end{equation}

The semi-norm (\ref{seminorm}) measures the energy or "smoothness" of the surface interpolant $s$, such that interpolants having a small semi-norm
are considered smoother than those having a large semi-norm. Following Duchon \cite{Duchon1977} and Green and Silverman \cite{Green1993}, the
solution to the problem (\ref{argmin}) is a function of the form
\begin{equation}\label{sform}
\widetilde s\left( {\boldsymbol{x}} \right) = {c_0} + {c_1}{\boldsymbol{x}} + \sum\limits_{i = 1}^{k \times N} {{\beta _i}K\left( {{{\left\| {{\boldsymbol{x}} - {{\boldsymbol{x}}_i}} \right\|}_2}} \right)},
\end{equation}
where $K$ is a real valued function defined on the kernel $K \in {\cal{K}}:{\mathbb{R}^{k \times N}} \times {\mathbb{R}^{k \times N}} \to
\mathbb{R}$, ${\left\| {\, \cdot \,} \right\|_2}$ is the Euclidian distance between the points $\boldsymbol x$ and $\boldsymbol x _i$, the
coefficients ${\beta_i} \in \mathbb{R}$ are constant real numbers and $\mathcal{P}\left( x \right)={c_0} + {c_1}{\boldsymbol{x}}$ is a global
polynomial function, usually considered of small degree. The points $\boldsymbol x _i$ are referred as centers of the Radial Basis Functions $K\left(
r \right) = {\cal{K}}\left( {\boldsymbol x,{\boldsymbol x_i}} \right)$, where the variable $r$ stands for ${{{\left\| {\boldsymbol x - {\boldsymbol
x_i}} \right\|}_2}}$.

In our approach we use the so called \textit{thinplate} kernel $K\left( r \right) = {r^2}\ln \left( {r + 1} \right)$. Ensuring that the interpolation
surface lies in the Beppo-Levi space $\widetilde s \in B{L_2}\left( {{\mathbb{R}^2}} \right)$ implies the following side conditions
\begin{equation}\label{sidecond1}
\sum\limits_{i = 1}^{k \times N} {{\beta_i}}  = \sum\limits_{i = 1}^{k \times N} {{x_i}{\beta_i}}  = \sum\limits_{i = 1}^{k \times N} {{y_i}{\beta_i}}  = 0
\end{equation}
and the constraints
\begin{equation}\label{sidecond2}
\sum\limits_{i = 1}^{k \times N} {{\beta_i}\mathcal{P}\left( {{\boldsymbol x_i}} \right)}  = 0.
\end{equation}
Considering that $\left\{ {{p_1},{p_2}} \right\}$ represents a basis for the polynomial $\mathcal{P}$  and $ \left\{ {{c_0},{c_1}} \right\}$ are the
coefficients that give the polynomial $\mathcal{P}\left( \boldsymbol x \right)$ in terms of this basis, the interpolation conditions
(\ref{interpcond}) with the side conditions (\ref{sidecond1}) and constraints (\ref{sidecond2}) lead to the following linear system to be solved for
the coefficients that specify the RBF
\begin{equation}\label{system}
\left( {\begin{array}{*{20}{c}}
K&P\\
{{P^T}}&0
\end{array}} \right)\left( {\begin{array}{*{20}{c}}
\beta \\
c
\end{array}} \right) = \left( {\begin{array}{*{20}{c}}
f\\
0
\end{array}} \right),
\end{equation}
where ${K_{ij}} = K\left( {{{\left\| {{\boldsymbol x_i} - {\boldsymbol x_j}} \right\|}_2}} \right),\; i,j = 1,...,k \times N$, ${P_{ij}} =
{p_j}\left( {{\boldsymbol x_i}} \right),\; i = 1,...,k \times N,\;j = 1,2$, $\beta = {\left( {{\beta_1},...,{\beta_{k \times N}}} \right)^T}$,  $c =
{\left( {{c_0},{c_1}} \right)^T}$, $f = {\left( {{f_1},...,{f_{k \times N}}} \right)^T}$. The zeros in (\ref{system}) denote matrices or vectors of
appropriate dimensions and 'T' stands for the transpose of a matrix or vector.
Solving the linear system (\ref{system}) determines the constant coefficients $\beta$ and the polynomial coefficients $c$ and hence the interpolant
surface $\widetilde s\left( \boldsymbol x \right)$.

The methodology presented herein leads to the following linear model (denoted in the following ARDMD-RBF model) for estimation of the flow field for
any time instance $t \in \left[ {{t_1},{t_N}} \right]$
\begin{equation}\label{lowmodel}
{v_{DMD}}^t\left( {\bf{x}} \right) = \sum\limits_{j = 1}^k {{b_j}^t{\phi _j}\left( {\bf{x}} \right)} ,\quad {b_j}^t = \widetilde s\left( {{\boldsymbol x_i}}  \right),\quad {{\boldsymbol x_i}}  \in \left\{ {1,...,k} \right\} \times \left[ {{t_1},{t_N}} \right],
\end{equation}
where ${{b_j}^t}$ are the interpolated coefficients, ${{\phi _j}\left( {\bf{x}} \right)}$ are the DMD basis functions, $k$ represents the number of
the DMD basis functions retained for the reduced order model and $t$ denotes any value of time in the interval $\left[ {{t_1},{t_N}} \right]$.

\section{Analysis of the reduced order model for Saint-Venant data}

\subsection{Acquisition of numerical data}

The test problem used in this paper is consisting of the nonlinear Saint Venant equations model (also called the shallow water equations
\cite{vre94}) in a channel on the rotating earth, associated with periodic boundary conditions in the  $\tilde x$-direction and solid wall boundary
condition in the $\tilde y$-direction:
\begin{equation}\label{sw1}
{\tilde u_{\tilde t}} + \tilde u{\tilde u_{\tilde x}} + \tilde v{\tilde u_{\tilde y}} + {\left( {g\tilde h} \right)_{\tilde x}} - \tilde f\tilde v = 0,
\end{equation}
\begin{equation}\label{sw2}
{\tilde v_{\tilde t}} + \tilde u{\tilde v_{\tilde x}} + \tilde v{\tilde v_{\tilde y}} + {\left( {g\tilde h} \right)_{\tilde y}} + \tilde f\tilde u = 0,
\end{equation}
\begin{equation}\label{sw3}
{\left( {g\tilde h} \right)_{\tilde t}} + {\left( {g\tilde h\tilde u} \right)_{\tilde x}} + {\left( {g\tilde h\tilde v} \right)_{\tilde y}} = 0,
\end{equation}
\begin{equation}\label{sw4}
\tilde u\left( {0,\tilde y,\tilde t} \right) = \tilde u\left( {{L_{\max }},\tilde y,\tilde t} \right),\;\tilde v\left( {\tilde x,0,\tilde t} \right) = \tilde v\left( {\tilde x,{D_{\max }},\tilde t} \right) = 0,
\end{equation}
where $\tilde u$ and $\tilde v$ are the velocity components in the $\tilde x$ and $\tilde y$ axis directions respectively, $g\tilde h$ is the
geopotential height, $\tilde h$ represents the depth of the fluid, $\tilde f$ is the Coriolis factor and $g$ is the acceleration of gravity. We
consider that the reference computational configuration is the rectangular $2D$ domain $\Omega = \left[ {0,{L_{\max }}} \right] \times \left[
{0,{D_{\max }}} \right]$. Subscripts represent the derivatives with respect to time and the streamwise and spanwise coordinate.

The Saint Venant equations have been used for a wide variety of hydrological and geophysical fluid dynamics phenomena such as tide-currents
\cite{Bryden2007}, pollutant dispersion \cite{Sportisse2000}, storm-surges or tsunami wave propagation \cite{kou88}. Early work on numerical methods
for solving the shallow water equations is described in Navon (1979) \cite{Navon1979}.

 We consider the model (\ref{sw1})-(\ref{sw4}) in a $\beta $-plane assumption detailed in \cite{Navon1983}, in which the effect of the Earth's sphericity is modeled by a linear variation in the Coriolis factor
\begin{equation}\label{coriolis}
\tilde f = {f_0} + \frac{\beta }{2}\left( {2\tilde y - {D_{\max }}} \right),
\end{equation}
where ${f_0},\beta $ are constants, ${L_{\max }},{D_{\max }}$ are the dimensions of the rectangular domain of integration.

The following initial condition introduced by Grammeltvedt \cite{Grammeltvedt1969} was adopted as the initial height field which propagates the
energy in wave number one, in the streamwise direction:
\begin{equation}\label{inicond1}
{h_0}\left( {\tilde x,\tilde y} \right) = {H_0} + {H_1}\tanh \left( {\frac{{9({D_{\max }}/2 - \tilde y)}}{{2{D_{\max }}}}} \right) + {H_2}\sin {\mkern 1mu} \left( {\frac{{2\pi \tilde x}}{{{L_{\max }}}}} \right){\cosh ^{ - 2}}\left( {\frac{{9({D_{\max }}/2 - \tilde y)}}{{{D_{\max }}}}} \right).
\end{equation}

Using the geostrophic relationship $\tilde u =  - {\tilde h_{\tilde y}}\left( {g/\tilde f} \right)$, $\tilde v = {\tilde h_{\tilde x}}\left(
{g/\tilde f} \right)$, the initial velocity fields are derived as:
\[{u_0}\left( {\tilde x,\tilde y} \right) =  - \frac{g}{{\tilde f}}\frac{{9{H_1}}}{{2{D_{\max }}}}\left( {{{\tanh }^2}\left( {\frac{{9{D_{\max }}/2 - 9\tilde y}}{{2{D_{\max }}}}} \right) - 1} \right) - \]
\begin{equation}\label{inicond2}
\frac{{18g}}{{\tilde f}}{H_2}\sinh \left( {\frac{{9{D_{\max }}/2 - 9\tilde y}}{{{D_{\max }}}}} \right)\frac{{\sin \left( {\frac{{2\pi \tilde x}}{{{L_{\max }}}}} \right)}}{{{D_{\max }}{{\cosh }^3}\left( {\frac{{9{D_{\max }}/2 - 9\tilde y}}{{{D_{\max }}}}} \right)}},
\end{equation}
\begin{equation}\label{inicond3}
{v_0}\left( {\tilde x,\tilde y} \right) = 2\pi {H_2}\frac{g}{{\tilde f{L_{\max }}}}\cos {\mkern 1mu} \left( {\frac{{2\pi \tilde x}}{{{L_{\max }}}}} \right){\cosh ^{ - 2}}\left( {\frac{{9({D_{\max }}/2 - \tilde y)}}{{{D_{\max }}}}} \right).
\end{equation}

The dimensional constants used for the above test model are
\[{f_0} = {10^{ - 4}}{s^{ - 1}},\quad \beta  = 1.5 \times {10^{ - 11}}{s^{ - 1}}{m^{ - 1}},\quad g = 10m{s^{ - 1}},\]
\begin{equation}\label{dimconst}
{{\rm{D}}_{\max }}{\rm{ = 44}} \times {\rm{1}}{{\rm{0}}^5}{\rm{m,}}\quad {{\rm{L}}_{\max }}{\rm{ = 6}} \times {\rm{1}}{{\rm{0}}^6}{\rm{m}},\quad {H_0} = 2 \times {\rm{1}}{{\rm{0}}^6}m,\quad {H_1} = 220m,\quad {H_2} = 133m.
\end{equation}

We have followed the approach used by Navon \cite{Navon1987}, which implements a two-stage finite-element Numerov-Galerkin method for integrating the
nonlinear shallow-water equations on a $\beta $-plane limited-area domain, for approximating the quadratic nonlinear terms that appear in the
equations of hydrological dynamics. This scheme when applied to meteorological and oceanographic problems gives an accurate phase propagation and
also handles nonlinearities well. The accuracy of temporal and spatial discretization equals $\mathcal{O}\left( {{k^2},{h^{p}}} \right)$, where $p$
varies in the interval $\left[ {4,8} \right]$. The training data comprises a number of $240$ unsteady solutions of the two-dimensional shallow water
equations model (\ref{sw1})-(\ref{inicond3}), at regularly spaced time intervals $\Delta t = 600s$ for each solution variable.

To measure the accuracy of the reduced shallow water model and to validate the numerical results with existing results in the literature, we
undertake a non-dimensional analysis of the shallow water model. Following \cite{Bare1996}, reference quantities of the dependent and independent
variables in the shallow water model are considered, i.e. the length scale ${L_{ref}} = {L_{\max }}$ and the reference units for the height and
velocity, respectively, are given by the initial conditions ${h_{ref}} = {h_0}$, ${u_{ref}} = {u_0}$. A typical time scale is also considered,
assuming the form ${t_{ref}} = {L_{ref}}/{u_{ref}}$. In order to make the system of equations (\ref{sw1})-(\ref{sw4}) non-dimensional, we define the
non-dimensional variables \[\left( {t,x,y} \right) = \left( {\tilde t/{t_{ref}},\tilde x/{L_{ref}},\tilde y/{L_{ref}}} \right),\quad \left( {h,u,v}
\right) = \left( {\tilde h/{h_{ref}},\tilde u/{u_{ref}},\tilde v/{u_{ref}}} \right).\] The numerical results are obtained and used in further
numerical experiments in dimensionless form.

\subsection{Advantages of ARDMD algorithm over classic approaches}

In this section, the efficiency of reduced order modelling based on the ARDMD algorithm is illustrated in comparison with  the classic DMD algorithm,
considering the evolution of the flow field along the integration time window. The first major advantage of the adaptive randomized DMD proposed in
this work is represented by the fact that the ARDMD algorithm produces a reduced order subspace of Ritz values, which has the same dimension as the
rank of randomized SVD function. In consequence, this procedure omits a further selection of the Ritz values associated with their DMD modes.

In case of the classic DMD algorithm, the superposition of all Koopman modes approximates the entire data sequence, but there are also modes that
have a weak contribution. The modes' selection, which is central in model reduction, constitutes the source of many discussions among modal
decomposition practitioners.  For instance, Jovanovic et al. \cite{Jovanovic2012} introduced a low-rank DMD algorithm to identify an a-priori
specified number of modes that provide optimal approximation of experimental or numerical snapshots at a certain time interval. Consequently, the
modes and frequencies that have strongest influence on the quality of approximation have been selected. An optimized DMD method was recently
introduced by Chen et al. \cite{Chen2012}, which tailors the decomposition to a desired number of modes. This method minimizes the total residual
over all data vectors and uses simulated annealing and quasi-Newton minimization iterative methods for selecting the frequencies.
 Tissot and his coworkers \cite{Tissot2014} propose a new energetic criterion for model reduction using Dynamic Mode Decomposition, that incorporates
 the growth rate of the modes.

%Ref. \cite{Bistrian2014} aimed to present a preliminary survey on DMD modes selection. In this previous investigation we introduced a novel selection
%method for the DMD modes, associated amplitudes and Ritz values to derive a DMD reduced-order model of the shallow water equations. We arranged the
%Koopman modes in descending order of the energy of the DMD modes weighted by the inverse of the Strouhal number. We eliminated the modes that
%contribute weakly to the data sequence based on the conservation of quadratic integral invariants by the reduced order flow. We address in this
%section the problem of identification of a small subset of DMD modes that yield an optimal truncated representation of the flow field in order to
%capture its most important dynamic structures. The recent investigation presented in \cite{bisnavon2016} has focused on the effects of modes
%selection in dynamic modes decomposition. We proposed an improved algorithm for selecting the dominant DMD modes from the flow field and we pointed
%out the problems for which the new criterion proved its efficiency.

\begin{figure}[h!]
\centering
\includegraphics[width=0.9\textwidth]{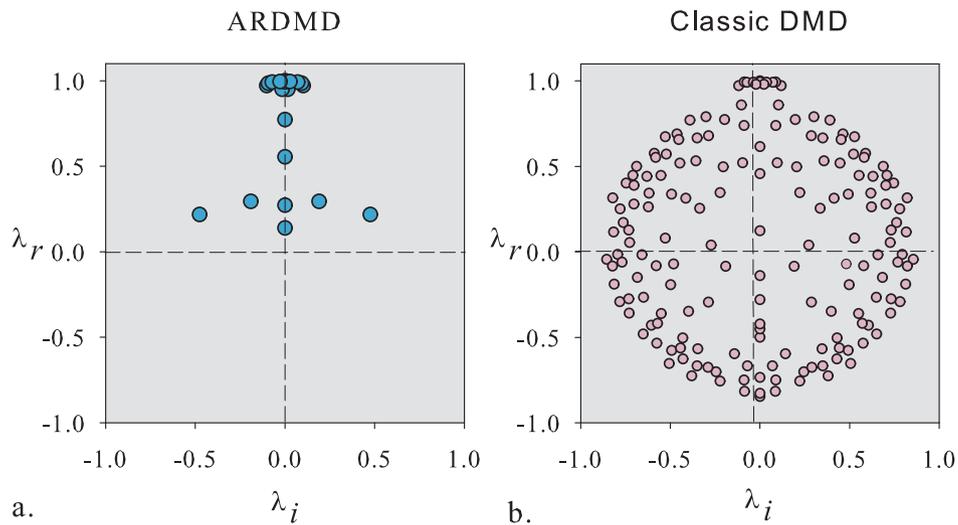}
\caption{The spectrum of DMD decomposition of geopotential height field $h$ in case of: a) ARDMD algorithm used in the present paper,  b) classic DMD algorithm used in \cite{bisnavon2016}.}\label{fig1}
\end{figure}

\begin{figure}[h!]
\centering
\includegraphics[width=0.9\textwidth]{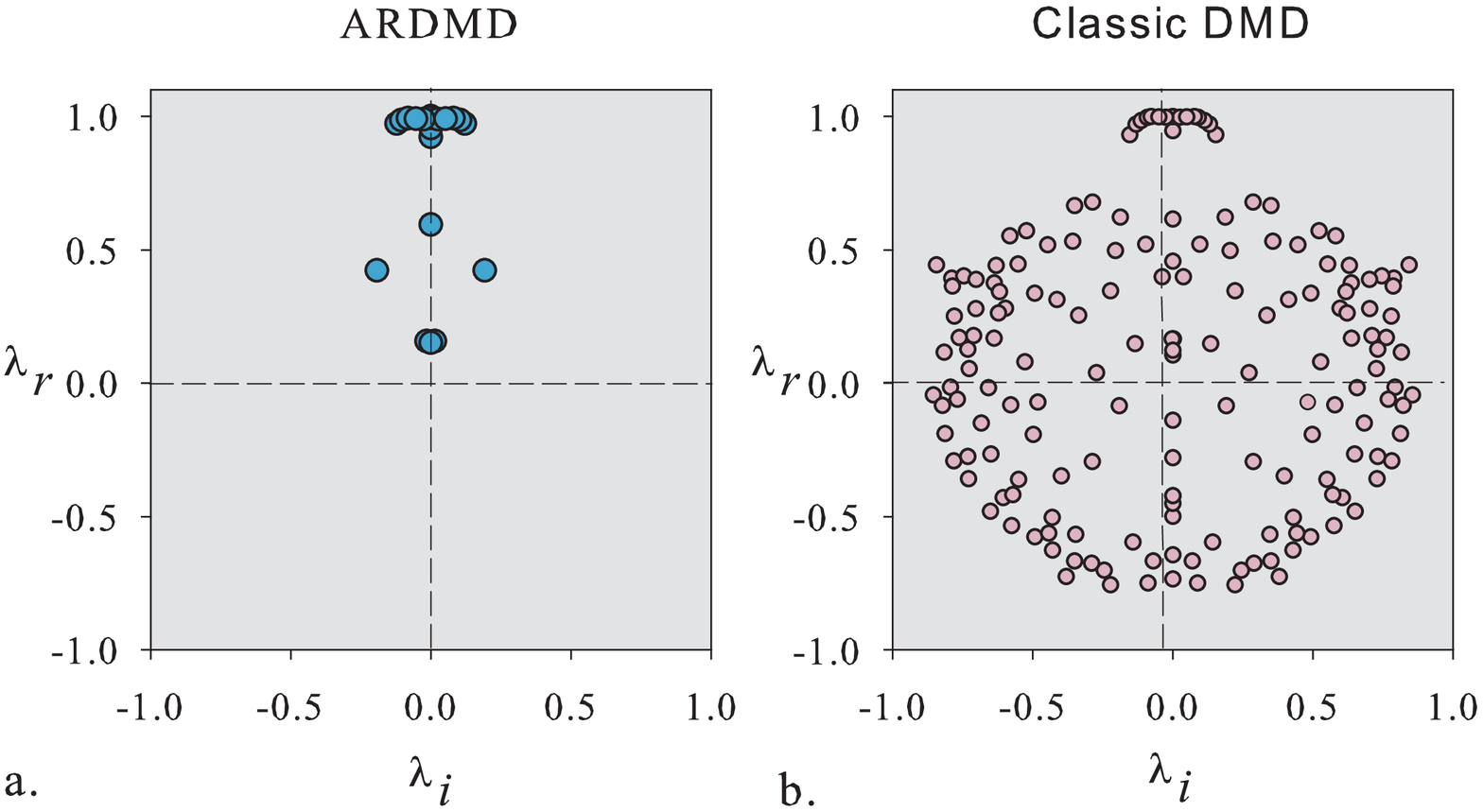}
\caption{The spectrum of DMD decomposition of streamwise field $u$ in case of: a) ARDMD algorithm used in the present paper,  b) classic DMD algorithm used in \cite{bisnavon2016}.}\label{fig2}
\end{figure}

\begin{figure}[h!]
\centering
\includegraphics[width=0.9\textwidth]{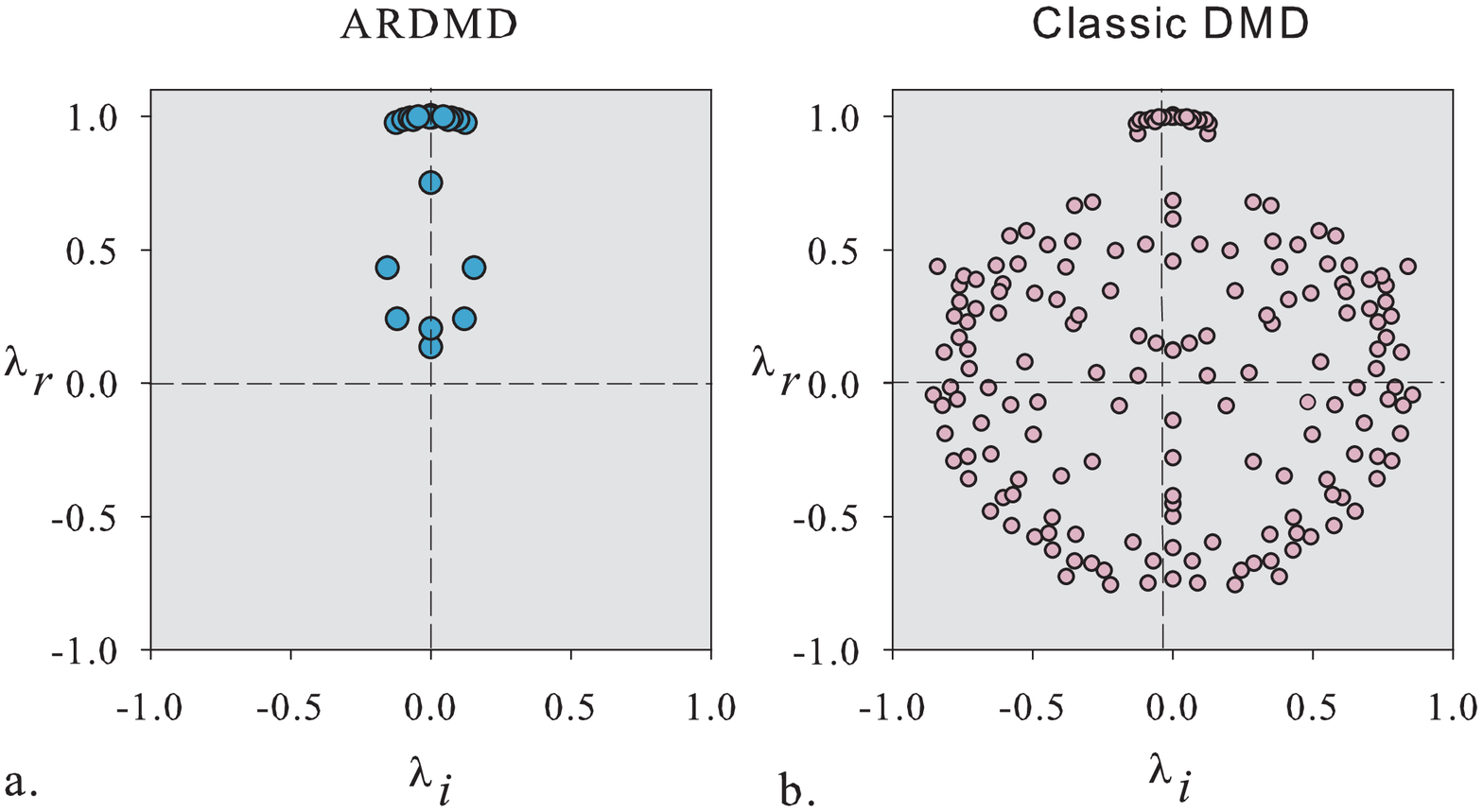}
\caption{The spectrum of DMD decomposition of spanwise field $v$ in case of: a) ARDMD algorithm used in the present paper,  b) classic DMD algorithm used in \cite{bisnavon2016}.}\label{fig3}
\end{figure}

To highlight the efficiency of the ARDMD method presented herein, we illustrate in Figures \ref{fig1}--\ref{fig3} the spectra of DMD decomposition of
geopotential height field $h$, streamwise field $u$ and spanwise field $v$, respectively, in case of the new ARDMD algorithm  and the classic DMD
algorithm  which was applied in \cite{bisnavon2016}.

\begin{figure}[h!]
\centering
\includegraphics[width=1\textwidth]{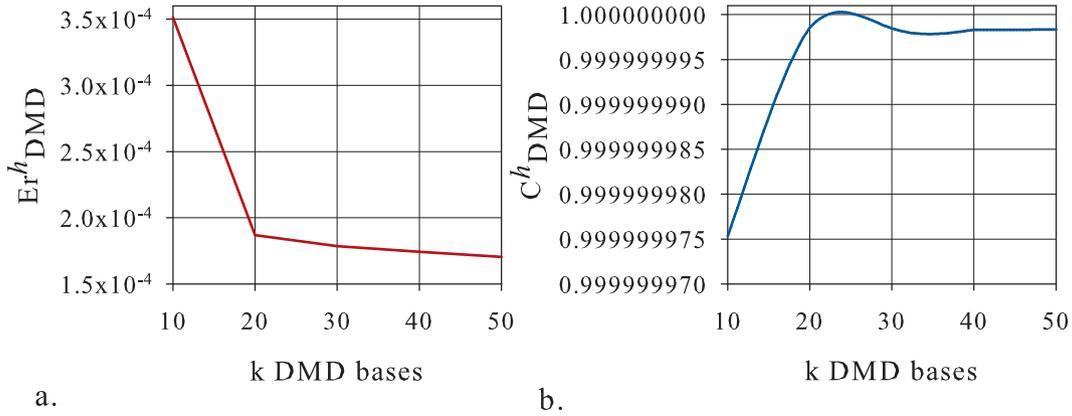}
\caption{ The process of evaluation of ROM target rank for $h$ field: a) The relative error of ARDMD computed as the retained number of dynamic modes,  b) The correlation coefficient of ARDMD computed as the retained number of dynamic modes.}\label{fig4}
\end{figure}

\begin{figure}[h!]
\centering
\includegraphics[width=1\textwidth]{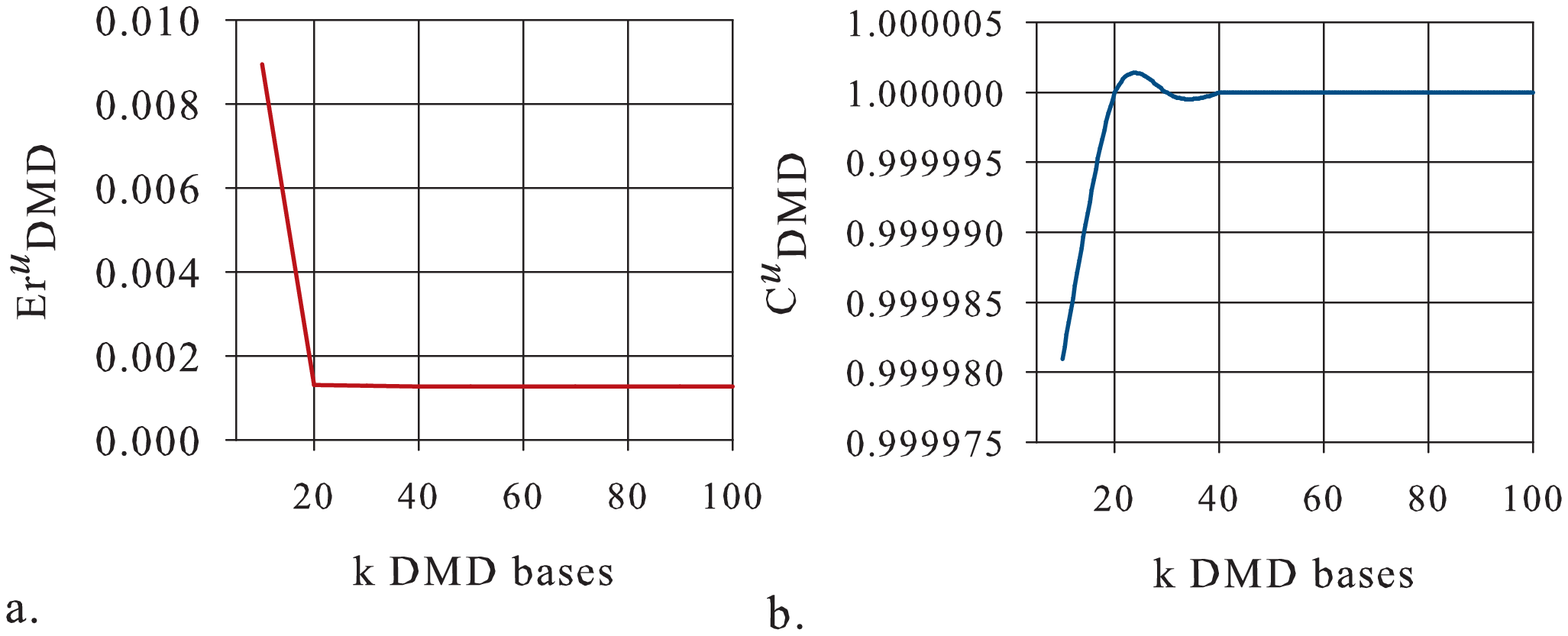}
\caption{ The process of evaluation of ROM target rank for $u$ field: a) The relative error of ARDMD computed as the retained number of dynamic modes,  b) The correlation coefficient of ARDMD computed as the retained number of dynamic modes.}\label{fig5}
\end{figure}

\begin{figure}[h!]
\centering
\includegraphics[width=1\textwidth]{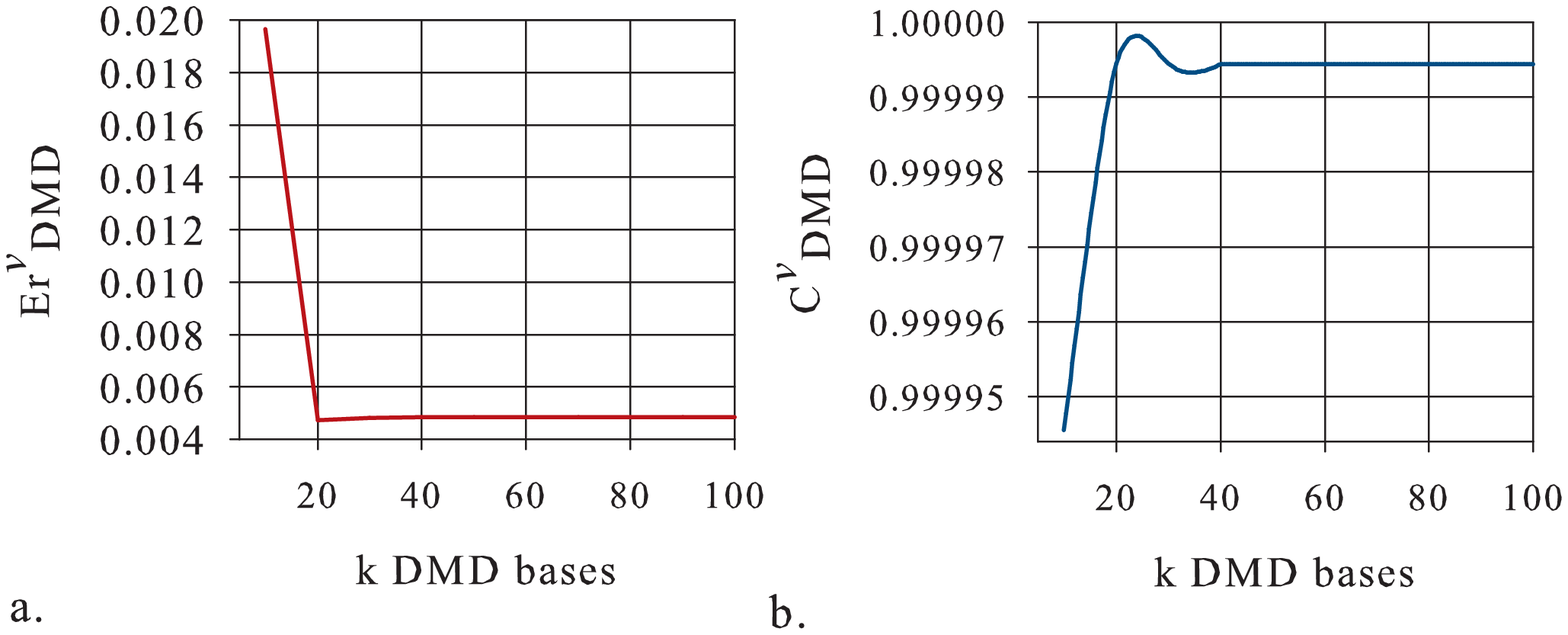}
\caption{ The process of evaluation of ROM target rank for $v$ field: a) The relative error of ARDMD computed as the retained number of dynamic modes,  b) The correlation coefficient of ARDMD computed as the retained number of dynamic modes.}\label{fig6}
\end{figure}

Obviously, when the classic DMD algorithm is applied, the practitioner has to address a modes' selection method like those mentioned above. Instead,
the randomized DMD algorithm (ARDMD) produces a significantly reduced size spectrum which elegantly incorporates the most influential modes. Figures
\ref{fig4}--\ref{fig6} present an insight of how ARDMD works.
%
%The randomized dynamic mode decomposition algorithm proposed in this paper solves the optimization problem (\ref{optimprob}) involving the method of
%sequential quadratic programming (SQP) \cite{Nocedal2006} which leads to the optimal low rank $k$ and associated DMD subspace $V_{DMD}$ where the
%most influential DMD bases live.
%
The optimization problem (\ref{optimprob}) is solved using a simulated annealing optimization routine which is detailed in \cite{SA}, based on
sequential quadratic programming (SQP) \cite{Nocedal2006}. This leads to the optimal low rank $k$ and associated DMD subspace $V_{DMD}$ where the
most influential DMD bases live. The rank of the reduced DMD model is automatically found such that the relative error of field reconstruction given
by Eq. (\ref{eror}) becomes sufficiently small and the correlation coefficient (\ref{corelation}) is sufficiently high. The optimal rank of the
reduced DMD model is the unique solution to the optimization problem (\ref{optimprob}). ARDMD algorithm produces subspaces of order $k=20$ selected
from $173$ DMD modes. A significant reduction of a factor of eight and a half is achieved for the representation of Saint-Venant fields $h$, $u$ and
$v$.

\begin{figure}[h!]
\centering
\includegraphics[width=1\textwidth]{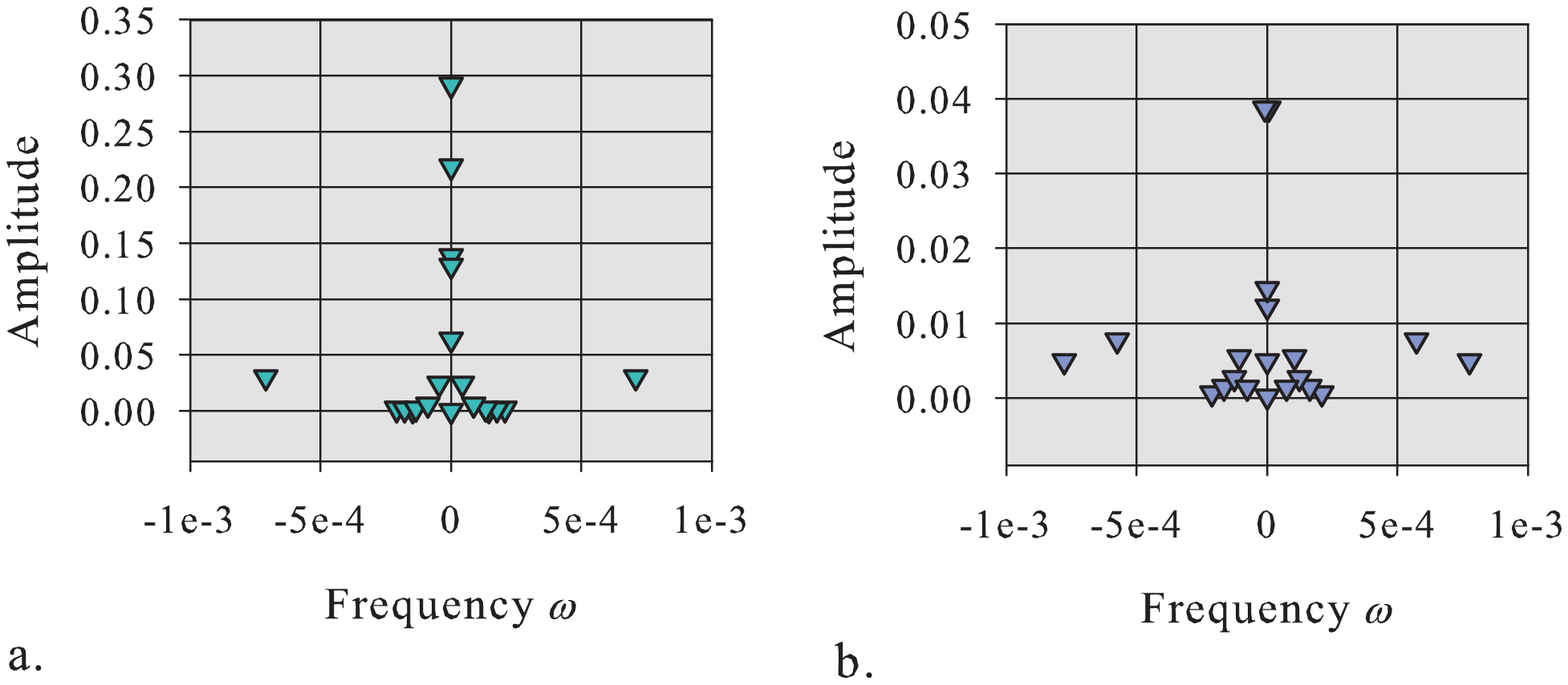}
\caption{ The amplitudes of DMD modes and associated frequencies obtained by dynamic mode decomposition of streamwise velocity field $u$ (a) and spanwise velocity field $v$ (b).}\label{fig7}
\end{figure}

The ARDMD algorithm presented herein is fully capable of determining the modal growth rates and the associated frequencies, which are illustrated in
Figure \ref{fig7} for velocity fields $u$, $v$, respectively. This is of major importance when is necessary to isolate modes with very high
amplitudes at lower frequencies or high frequency modes having lower amplitudes.

The numerical results provided by the ARDMD algorithm are presented in Table \ref{table1}.
\begin{table}[h!]
\caption{The relative error $E{r}_{DMD}$, the correlation coefficient $C_{DMD}$ and the reduced order modelling rank $k$
 obtained from ARDMD modal decomposition.}
\centering
\begin{threeparttable}
\tabsize
\begin{tabular}{cccc}
\toprule
\begin{small}  Flow   \end{small}& \begin{small} The relative  \end{small}  & \begin{small} The correlation \end{small} & \begin{small} The ROM \end{small}\\
\begin{small}   field   \end{small}& \begin{small} error   \end{small}  & \begin{small} coefficient \end{small} & \begin{small} rank  \end{small}\\
\midrule
 $h\left( {x,y} \right)$ & $E{r^h}_{DMD} = 1.7909 \times {10^{ - 4}} $ & ${C^h}_{DMD} = 0.99999$ & $k=20$\\
 $u\left( {x,y} \right)$ & $E{r^u}_{DMD} = 1.4738 \times {10^{ - 3}}$ & ${C^u}_{DMD} = 0.99999$ & $k=20$\\
 $v\left( {x,y} \right)$ & $E{r^v}_{DMD} = 4.5316 \times {10^{ - 3}}$ & ${C^v}_{DMD} = 0.99999$ & $k=20$\\
\bottomrule
\end{tabular} %\label{table1}
\begin{tablenotes}
      \small
      \item ARDMD, adaptive randomized dynamic mode decompositon, ROM, reduced order model.
    \end{tablenotes}
\end{threeparttable}\label{table1}
\end{table}

A comparison of the reduced order modelling rank, in the case of several DMD based modal decomposition methods associated with certain modes'
selection criteria proposed in our previous investigations and novel ARDMD technique is presented in Table \ref{table2}. Ref. \cite{Bistrian2014}
aimed to present a preliminary survey on DMD modes selection. We proposed a framework for modal decomposition of $2D$ flows, when numerical or
experimental data are captured with large time steps. Key innovations for the DMD-based ROM introduced in \cite{Bistrian2014} are the use of the
Moore–-Penrose pseudoinverse in the DMD computation that produced an accurate result and a novel selection method for the DMD modes. Unlike the
classic algorithm, we arrange the Koopman modes in descending order of the energy of the DMD modes weighted by the inverse of the Strouhal number. We
eliminate the modes that contribute weakly to the data sequence based on the conservation of quadratic integral invariants \cite{Navon1986} by the
reduced order flow. The resulting optimization problem was solved by employing sequential quadratic programming (SQP) \cite{Nocedal2006}.

In \cite{aleksey2015} we proposed a new framework for dynamic mode decomposition based on the reduced Schmid operator. We investigated a variant of
DMD algorithm and we explored the selection of the modes based on sorting them in decreasing order of their amplitudes.
 This procedure works well for models without modes that are very rapidly damped, having very high amplitudes. Therefore the selection of modes based
 on their amplitude is effective only in certain situations, as reported also by Noack et al. \cite{Noak2011}.

The investigation recently presented in \cite{bisnavon2016} has focused on the effects of modes selection in dynamic mode decomposition. We proposed
a new vector filtering criterion for dynamic modes selection that is able to extract dynamically relevant flow features of time-resolved experimental
or numerical data. The algorithm related in \cite{bisnavon2016} proposed a dynamic filtering criterion for which the amplitude of any mode is
weighted by its growth rate. This method proved to be perfectly adapted to the flow dynamics, in identification of the most influential modes.

\begin{table}[h!]
\caption{The reduced order modelling rank $k$ and the relative error $E{r^h}_{DMD}$ in case of several DMD based modal decomposition methods.}
\centering
\begin{threeparttable}
\tabsize
\begin{tabular}{cccc}
\toprule
\begin{small}  energetic   \end{small}& \begin{small} reduced Schmid  \end{small}  & \begin{small} dynamic \end{small} & \begin{small} ARDMD algorithm \end{small}\\
\begin{small}   DMD (\cite{Bistrian2014})  \end{small}& \begin{small} operator DMD  (\cite{aleksey2015}) \end{small}  & \begin{small} DMD (\cite{bisnavon2016}) \end{small} & \begin{small} the present research \end{small}\\
\midrule
 $13$ & $19$ & $11$ & $20$\\
 $1.19\times 10^{-3}$ & $2.683\times 10^{-4}$ & $2.5785\times 10^{-4}$ & $1.7909\times 10^{-4} $\\
\bottomrule
\end{tabular} %\label{table1}
\begin{tablenotes}
      \small
      \item DMD, dynamic mode decomposition, ARDMD, adaptive randomized dynamic mode decompositon.
    \end{tablenotes}
\end{threeparttable}\label{table2}
\end{table}

Data presented in Table \ref{table2} argues the efficiency of the novel ARDMD method. Although the previous techniques detailed in
\cite{Bistrian2014, aleksey2015, bisnavon2016} lead to a reduced number of retained modes, there are still missing modes that would contribute to
data approximation. Hence the relative error of flow reconstruction by the reduced order model is the best in the case of randomized dynamic mode
decomposition. Producing a slightly larger model rank $k$ than the previous algorithms, the great advantage of adaptive randomized DMD (ARDMD) is
that omits the efforts of implementing an additional criterion of influential modes'selection, they being selected automatically.

\begin{figure}[h!]
\centering
\includegraphics[width=0.8\textwidth]{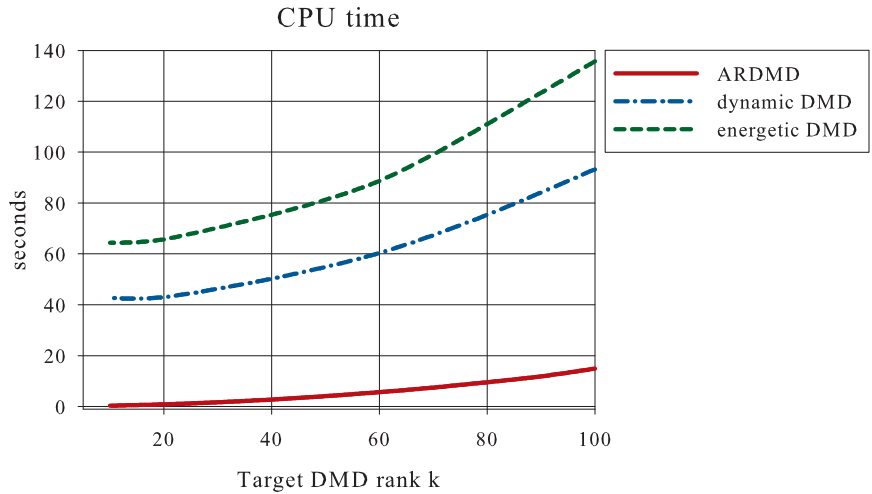}
\caption{ CPU time gained by applying several DMD techniques: ARDMD, adaptive randomized dynamic mode decomposition (present research), dynamic DMD method presented  in \cite{bisnavon2016}, energetic DMD method employed in \cite{Bistrian2014}.}\label{fig8}
\end{figure}

Thus a significant reduction in computational time is also achieved. The CPU time gained by applying several DMD techniques is presented in Figure
\ref{fig8}. By employing the ARDMD algorithm in comparison with dynamic DMD \cite{bisnavon2016} and energetic DMD \cite{Bistrian2014}, the
computational complexity of the reduced order model is diminished from the very beginning by two times and three times, respectively, as illustrated
in Figure \ref{fig8}.

\begin{figure}[h!]
\centering
\includegraphics[width=0.9\textwidth]{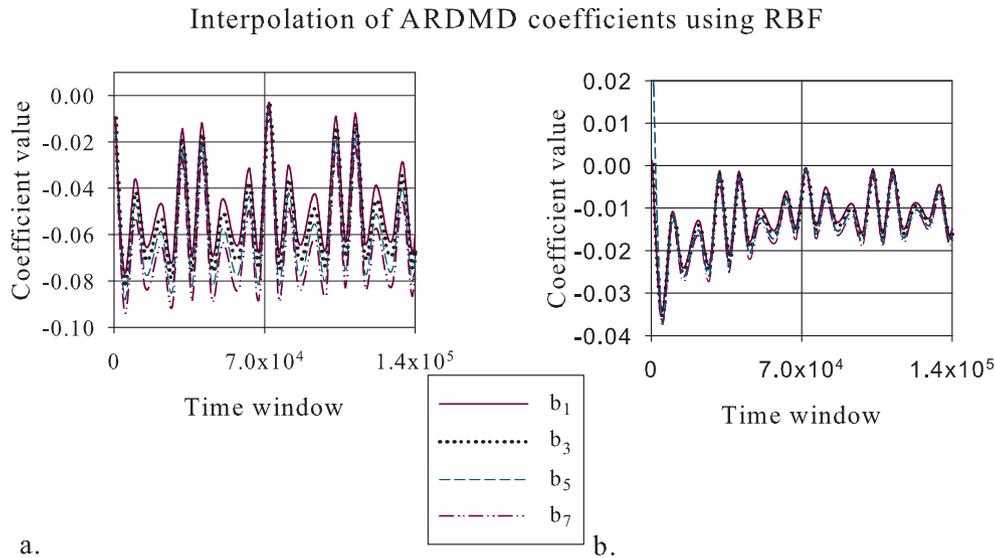}
\caption{The coefficients of the reduced order models (a) ${h_{DMD}}\left( {x,y} \right)$ and (b) ${u_{DMD}}\left( {x,y} \right)$ obtained by radial basis functions (RBF) interpolation, ${b_j} = {a_j}\lambda _j^{t - 1}$ for $j = 1,...,k$.}\label{fig9}
\end{figure}

Estimation of low order model coefficients by interpolation, in cases of non-intrusive data, represent a cost effective solution, as has been also
reported in the literature by Raisee et al. \cite{Raisee2015}, Peherstorfer and Willcox \cite{Peher2016}, Lin et al. \cite{Lin2016}. The coefficients
of the reduced order models of state solutions $\left( {{h_{DMD}},{u_{DMD}},{v_{DMD}}} \right)\left( {x,y} \right)$ have been estimated for entire
time window by interpolating the DMD computed coefficients using radial basis functions (RBF) discussed in Section \ref{RBF}. They are depicted in
Figure \ref{fig9}.

\begin{figure}[h!]
\centering
\includegraphics[width=0.7\textwidth]{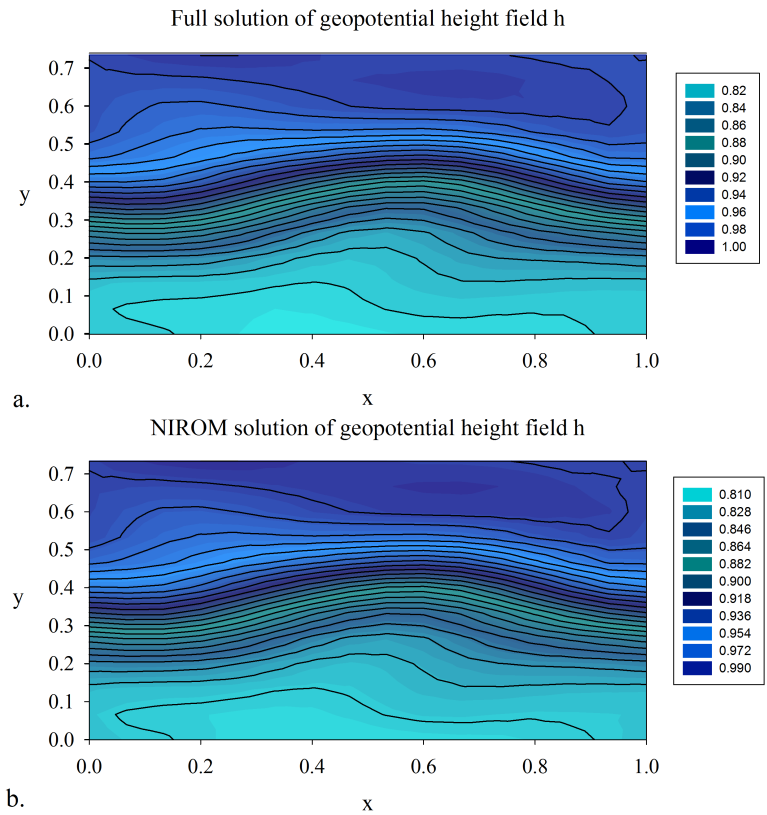}
\caption{a) Full solution of geopotential height field; b) NIROM solution of geopotential height field.}\label{fig10}
\end{figure}

\begin{figure}[h!]
\centering
\includegraphics[width=0.75\textwidth]{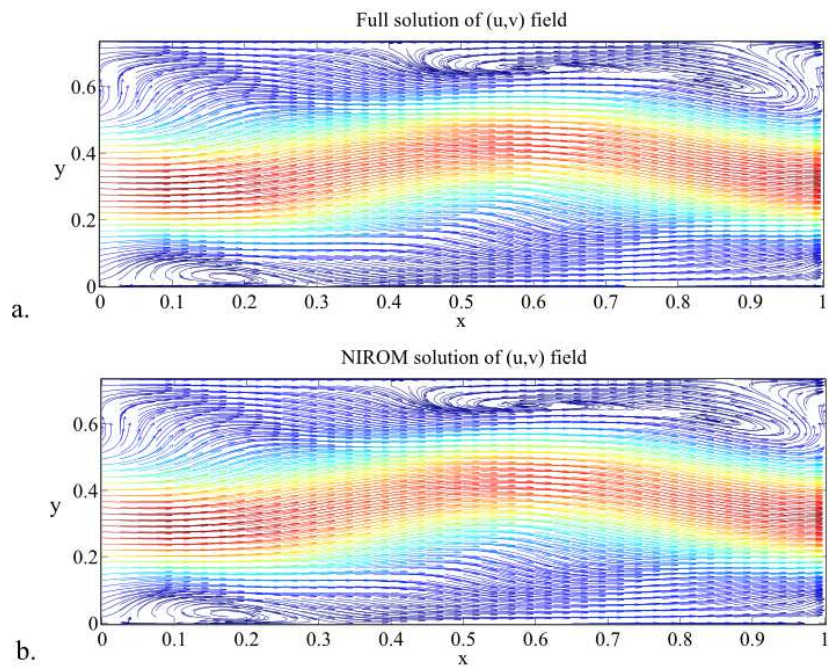}
\caption{a) Full solution of $\left( {u,v} \right)$ field; b) NIROM solution of $\left( {u,v} \right)$ field.}\label{fig11}
\end{figure}

Using the ARDMD algorithm, we obtain the non-intrusive reduced order models (NIROMs) of the state solutions $\left( {h,u,v} \right)\left( {x,y,t}
\right)$. The validity of the methodology introduced in this paper is checked by comparing how the NIROMs assess the full solution fields given by
the experimental data at time instance $181$, in Figures \ref{fig10} and \ref{fig11}, respectively. We applied a normalization condition such that
the maximum amplitude of the physical components $\left( {h,u,v} \right)\left( {x,y,t} \right)$ fields over the $\left( {x,y} \right)$ stations is
unity. The NIROM models employ the Ritz values represented in Figures \ref{fig1}-\ref{fig3} and associated DMD modes and amplitudes determined by the
ARDMD algorithm.

\begin{figure}[h!]
\centering
\includegraphics[width=0.75\textwidth]{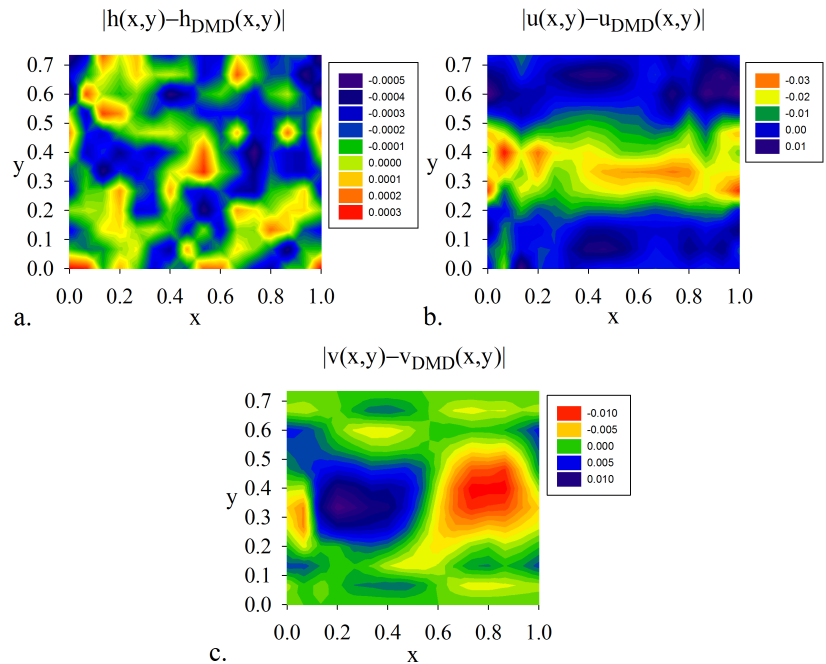}
\caption{Local error between full SWE solution and NIROM of geopotential height field (a), NIROM of streamwise velocity field (b), NIROM of spanwise velocity field (c), respectively.}\label{fig12}
\end{figure}

The local error between the full SWE solution and NIROM solution, respectively, at time instance $181$ is provided in Figure \ref{fig12}.

\begin{figure}[h!]
\centering
\includegraphics[width=0.8\textwidth]{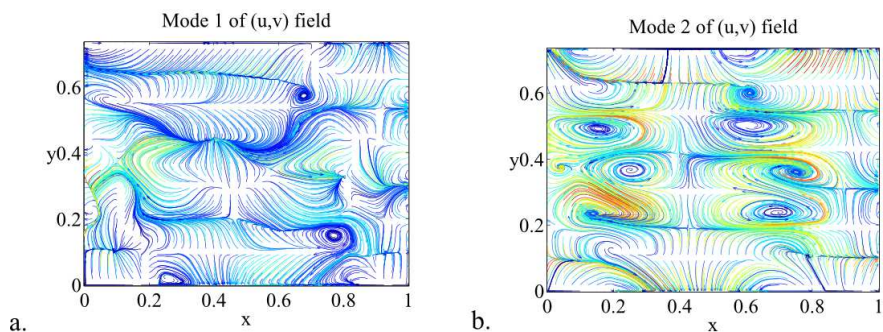}
\caption{The first two DMD modes of the $\left( {u,v} \right)$ field.}\label{fig13}
\end{figure}

\begin{figure}[h!]
\centering
\includegraphics[width=0.8\textwidth]{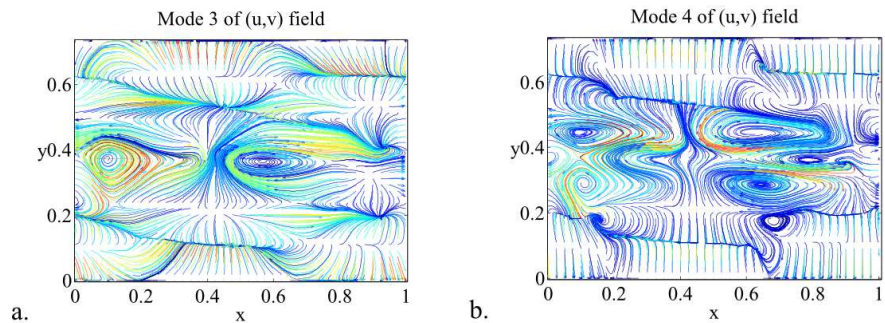}
\caption{The third and fourth DMD modes of the $\left( {u,v} \right)$ field.}\label{fig14}
\end{figure}

The coherent structures in the $\left( {u,v} \right)$ field can be visualized as local vortices in the first DMD modes, which are illustrated in
Figure \ref{fig13} and Figure \ref{fig14}.

The flow reconstructions by non-intrusive reduced order models (NIROMs) presented in Figure \ref{fig10} and Figure \ref{fig11} are very close to
experimental data snapshots, comparing the solution of SWE flow field after $181$ steps. The values of the correlation coefficients provided in Table
\ref{table1} greater than $99\%$ confirm the validity of the NIROMs. The similarity between the characteristics of the flow field and those obtained
by the NIROMs validates the method presented here and certifies that the improved ARDMD method can be applied successfully to model reduction of $2D$
flows.

For the problem investigated here, the novel ARDMD method shows satisfactory performances and provides a higher degree of accuracy for flow linear
reduced order model.

\section{Summary and Conclusions}

In this paper we have proposed a framework for reduced order modelling of non-intrusive data with application to $2D$ flows. To overcome the
inconveniences of intrusive model order reduction usually derived by combining the POD and the Galerkin projection methods, we developed a novel
technique based on randomized dynamic mode decomposition as a fast and accurate option in model order reduction of non-intrusive data originating
from Saint-Venant systems. We derived a non-intrusive approach to obtain a reduced order linear model of the flow dynamics based on dynamic mode
decomposition of experimental data in association with the efficient radial basis function interpolation technique.

To the best of our knowledge, the present paper is the first work that introduced the randomized dynamic mode decomposition algorithm with
application to fluid dynamics, after the randomized SVD algorithm recently introduced by Erichson and Donovan \cite{Erichson2016} for processing of
high resolution videos.

Several key innovations have been introduced in the present paper:
\begin{itemize}
\item We endow the dynamic mode decomposition (DMD) algorithm with a randomized singular value decomposition (RSVD) algorithm.
\item We gain a fast and accurate adaptive randomized DMD algorithm (ARDMD), exploiting an efficient low rank DMD model of input data.
\item The rank of the reduced DMD model is given as the unique solution of an optimization problem whose constraints are a sufficiently small
    relative error of data reconstruction and a sufficiently high correlation coefficient between the experimental data and the DMD solution.
\end{itemize}

The major advantages of the adaptive randomized dynamic mode decomposition (ARDMD) proposed in this work are:
\begin{itemize}
\item This method provides an efficient tool in developing the linear model of a complex flow field described by non-intrusive (or experimental)
    data.
\item This method does not require an additional selection algorithm of the DMD modes. ARDMD produces a reduced order subspace of Ritz values,
    having the same dimension as the rank of randomized SVD function, where the most influential DMD modes live.

\item We gain a significantly reduction of CPU time in computation of ROMs for massive numerical data.

\item Combining the randomized dynamic mode decomposition with radial basis function (RBF) interpolation we have derived a reduced order model for
    estimating the flow behavior in the real time window. The non-intrusive reduced order model (NIROM) presents satisfactory performances in flow
    reconstruction.

\item Analyzing the modal growth rates and the associated frequencies is an instance of apprehending the flow dynamics. This is of major importance
    when is necessary to isolate modes with very high amplitudes at lower frequencies or high frequency modes having lower amplitudes. Thus this
    paper outlines steps towards hydrodynamic stability analysis and flow control with potential applications.

\end{itemize}

To highlight the performances of the proposed methodology we performed a comparison of the reduced order modelling rank, in the case of several DMD
based modal decomposition methods associated with certain modes’ selection criteria and novel ARDMD technique presened herein. The numerical results
argued the efficiency of the novel ARDMD method.

We emphasized the excellent behavior of the NIROMs developed in this paper by comparing the computed shallow water solution with the experimentally
measured profiles and we found a close agreement. In addition, we performed a qualitative analysis of the reduced order models by  correlation
coefficients and local errors.

There are a number of interesting directions that arise from this work. First, it will be a natural extension to apply the proposed algorithm to
high-dimensional data originating from fluid dynamics and  oceanographic/atmospheric measurements. The methodology presented here offers the main
advantage of deriving a reduced order model capable to provide a variety of information describing the behavior of the flow field. A future extension
of this research will address an efficient numerical approach for modal decomposition of swirling flows, where the full mathematical model implies
more sophisticated relations at domain boundaries that must be satisfied by the reduced order  model also.

%%%%%%%%%%%%%%%%%%%%%%%%%

%\ack This class file was developed by Sunrise Setting Ltd,
%Brixham, Devon, UK. Website:\\
%\href{http://www.sunrise-setting.co.uk}{\texttt{www.sunrise-setting.co.uk}}

\bibliographystyle{wileyj}
%\bibliography{DiaJSC}
%\IfFileExists{\jobname.bbl}{}
% {\typeout{}
%  \typeout{******************************************}
%  \typeout{** Please run "bibtex \jobname" to optain}
% \typeout{** the bibliography and then re-run LaTeX}
% \typeout{** twice to fix the references!}
% \typeout{******************************************}
% \typeout{}
%  }
%%

\end{document}